\documentclass[]{amsart}
\usepackage{amsthm}
\usepackage{amsfonts}

\catcode`\@=11\relax
\newwrite\@unused
\def\typeout#1{{\let\protect\string\immediate\write\@unused{#1}}}
\def\@nnil{\@nil}
\def\@empty{}
\def\@psdonoop#1\@@#2#3{}
\def\@psdo#1:=#2\do#3{\edef\@psdotmp{#2}\ifx\@psdotmp\@empty \else
    \expandafter\@psdoloop#2,\@nil,\@nil\@@#1{#3}\fi}
\def\@psdoloop#1,#2,#3\@@#4#5{\def#4{#1}\ifx #4\@nnil \else
       #5\def#4{#2}\ifx #4\@nnil \else#5\@ipsdoloop #3\@@#4{#5}\fi\fi}
\def\@ipsdoloop#1,#2\@@#3#4{\def#3{#1}\ifx #3\@nnil 
       \let\@nextwhile=\@psdonoop \else
      #4\relax\let\@nextwhile=\@ipsdoloop\fi\@nextwhile#2\@@#3{#4}}
\def\@tpsdo#1:=#2\do#3{\xdef\@psdotmp{#2}\ifx\@psdotmp\@empty \else
    \@tpsdoloop#2\@nil\@nil\@@#1{#3}\fi}
\def\@tpsdoloop#1#2\@@#3#4{\def#3{#1}\ifx #3\@nnil 
       \let\@nextwhile=\@psdonoop \else
      #4\relax\let\@nextwhile=\@tpsdoloop\fi\@nextwhile#2\@@#3{#4}}
\def\psdraft{
	\def\@psdraft{0}
}
\def\psfull{
	\def\@psdraft{100}
}
\psfull
\newif\if@prologfile
\newif\if@postlogfile
\newif\if@bbllx
\newif\if@bblly
\newif\if@bburx
\newif\if@bbury
\newif\if@height
\newif\if@width
\newif\if@rheight
\newif\if@rwidth
\newif\if@clip
\def\@p@@sclip#1{\@cliptrue}
\def\@p@@sfile#1{
		   \def\@p@sfile{#1}
}
\def\@p@@sfigure#1{\def\@p@sfile{#1}}
\def\@p@@sbbllx#1{
		\@bbllxtrue
		\dimen100=#1
		\edef\@p@sbbllx{\number\dimen100}
}
\def\@p@@sbblly#1{
		\@bbllytrue
		\dimen100=#1
		\edef\@p@sbblly{\number\dimen100}
}
\def\@p@@sbburx#1{
		\@bburxtrue
		\dimen100=#1
		\edef\@p@sbburx{\number\dimen100}
}
\def\@p@@sbbury#1{
		\@bburytrue
		\dimen100=#1
		\edef\@p@sbbury{\number\dimen100}
}
\def\@p@@sheight#1{
		\@heighttrue
		\dimen100=#1
   		\edef\@p@sheight{\number\dimen100}
}
\def\@p@@swidth#1{
		\@widthtrue
		\dimen100=#1
		\edef\@p@swidth{\number\dimen100}
}
\def\@p@@srheight#1{
		\@rheighttrue
		\dimen100=#1
		\edef\@p@srheight{\number\dimen100}
}
\def\@p@@srwidth#1{
		\@rwidthtrue
		\dimen100=#1
		\edef\@p@srwidth{\number\dimen100}
}
\def\@p@@sprolog#1{\@prologfiletrue\def\@prologfileval{#1}}
\def\@p@@spostlog#1{\@postlogfiletrue\def\@postlogfileval{#1}}
\def\@cs@name#1{\csname #1\endcsname}
\def\@setparms#1=#2,{\@cs@name{@p@@s#1}{#2}}
\def\ps@init@parms{
		\@bbllxfalse \@bbllyfalse
		\@bburxfalse \@bburyfalse
		\@heightfalse \@widthfalse
		\@rheightfalse \@rwidthfalse
		\def\@p@sbbllx{}\def\@p@sbblly{}
		\def\@p@sbburx{}\def\@p@sbbury{}
		\def\@p@sheight{}\def\@p@swidth{}
		\def\@p@srheight{}\def\@p@srwidth{}
		\def\@p@sfile{}
		\def\@p@scost{10}
		\def\@sc{}
		\@prologfilefalse
		\@postlogfilefalse
		\@clipfalse
}
\def\parse@ps@parms#1{
	 	\@psdo\@psfiga:=#1\do
		   {\expandafter\@setparms\@psfiga,}}
\newif\ifno@bb
\newif\ifnot@eof
\newread\ps@stream
\def\bb@missing{
	\typeout{psfig: searching \@p@sfile \space  for bounding box}
	\openin\ps@stream=\@p@sfile
	\no@bbtrue
	\not@eoftrue
	\catcode`\%=12
	\loop
		\read\ps@stream to \line@in
		\global\toks200=\expandafter{\line@in}
		\ifeof\ps@stream \not@eoffalse \fi
		\@bbtest{\toks200}
		\if@bbmatch\not@eoffalse\expandafter\bb@cull\the\toks200\fi
	\ifnot@eof \repeat
	\catcode`\%=14
}	
\catcode`\%=12
\newif\if@bbmatch
\def\@bbtest#1{\expandafter\@a@\the#1
\long\def\@a@#1
\long\def\bb@cull#1 #2 #3 #4 #5 {
	\dimen100=#2 bp\edef\@p@sbbllx{\number\dimen100}
	\dimen100=#3 bp\edef\@p@sbblly{\number\dimen100}
	\dimen100=#4 bp\edef\@p@sbburx{\number\dimen100}
	\dimen100=#5 bp\edef\@p@sbbury{\number\dimen100}
	\no@bbfalse
}
\catcode`\%=14
\def\compute@bb{
		\no@bbfalse
		\if@bbllx \else \no@bbtrue \fi
		\if@bblly \else \no@bbtrue \fi
		\if@bburx \else \no@bbtrue \fi
		\if@bbury \else \no@bbtrue \fi
		\ifno@bb \bb@missing \fi
		\ifno@bb \typeout{FATAL ERROR: no bb supplied or found}
			\no-bb-error
		\fi
		\count203=\@p@sbburx
		\count204=\@p@sbbury
		\advance\count203 by -\@p@sbbllx
		\advance\count204 by -\@p@sbblly
		\edef\@bbw{\number\count203}
		\edef\@bbh{\number\count204}
}
\def\in@hundreds#1#2#3{\count240=#2 \count241=#3
		     \count100=\count240	
		     \divide\count100 by \count241
		     \count101=\count100
		     \multiply\count101 by \count241
		     \advance\count240 by -\count101
		     \multiply\count240 by 10
		     \count101=\count240	
		     \divide\count101 by \count241
		     \count102=\count101
		     \multiply\count102 by \count241
		     \advance\count240 by -\count102
		     \multiply\count240 by 10
		     \count102=\count240	
		     \divide\count102 by \count241
		     \count200=#1\count205=0
		     \count201=\count200
			\multiply\count201 by \count100
		 	\advance\count205 by \count201
		     \count201=\count200
			\divide\count201 by 10
			\multiply\count201 by \count101
			\advance\count205 by \count201
		     \count201=\count200
			\divide\count201 by 100
			\multiply\count201 by \count102
			\advance\count205 by \count201
		     \edef\@result{\number\count205}
}
\def\compute@wfromh{
		\in@hundreds{\@p@sheight}{\@bbw}{\@bbh}
		\edef\@p@swidth{\@result}
}
\def\compute@hfromw{
		\in@hundreds{\@p@swidth}{\@bbh}{\@bbw}
		\edef\@p@sheight{\@result}
}
\def\compute@handw{
		\if@height 
			\if@width
			\else
				\compute@wfromh
			\fi
		\else 
			\if@width
				\compute@hfromw
			\else
				\edef\@p@sheight{\@bbh}
				\edef\@p@swidth{\@bbw}
			\fi
		\fi
}
\def\compute@resv{
		\if@rheight \else \edef\@p@srheight{\@p@sheight} \fi
		\if@rwidth \else \edef\@p@srwidth{\@p@swidth} \fi
}
\def\compute@sizes{
	\compute@bb
	\compute@handw
	\compute@resv
}
\def\psfig#1{\vbox {
	\ps@init@parms
	\parse@ps@parms{#1}
	\compute@sizes
	\ifnum\@p@scost<\@psdraft{
		\typeout{psfig: including \@p@sfile \space }
		\special{ps::[begin] 	\@p@swidth \space \@p@sheight \space
				\@p@sbbllx \space \@p@sbblly \space
				\@p@sbburx \space \@p@sbbury \space
				startTexFig \space }
		\if@clip{
			\typeout{(clip)}
			\special{ps:: \@p@sbbllx \space \@p@sbblly \space
				\@p@sbburx \space \@p@sbbury \space
				doclip \space }
		}\fi
		\if@prologfile
		    \special{ps: plotfile \@prologfileval \space } \fi
		\special{ps: plotfile \@p@sfile \space }
		\if@postlogfile
		    \special{ps: plotfile \@postlogfileval \space } \fi
		\special{ps::[end] endTexFig \space }
		\vbox to \@p@srheight true sp{
			\hbox to \@p@srwidth true sp{
				\hfil
			}
		\vfil
		}
	}\else{
		\vbox to \@p@srheight true sp{
		\vss
			\hbox to \@p@srwidth true sp{
				\hss
				\@p@sfile
				\hss
			}
		\vss
		}
	}\fi
}}
\catcode`\@=12\relax

\def\func#1{\mathop{\rm #1}}

\newtheorem{theorem}{Theorem}[section]

\newtheorem{lemma}{Lemma}[section]
\newtheorem{definition}{Definition}[section]

\newtheorem{remark}{Remark}[section]

\begin{document}

\title[Infinite Periodic Discrete Minimal Surfaces]{Infinite Periodic 
Discrete Minimal Surfaces Without Self-Intersections}
\author{Wayne Rossman}
\begin{abstract}
A triangulated piecewise-linear minimal surface in Euclidean $3$-space 
${\mathbb{R}^{3}}$ defined 
using a variational characterization is critical for area 
amongst all continuous piecewise-linear variations with compact support 
that preserve the simplicial structure.  We explicitly construct examples of 
such surfaces that are embedded and are periodic in three independent 
directions of ${\mathbb{R}^{3}}$.  
\end{abstract}

\maketitle

\begin{center}
Department of Mathematics, Faculty of Science, \\ 
Kobe University, Rokko, Kobe 657-8501, Japan \\ 
wayne@math.kobe-u.ac.jp \\ 
http://math.kobe-u.ac.jp/HOME/wayne/wayne.html \\
Tel.: +81-78-803-5623, FAX: +81-78-803-5610 
\end{center}

\begin{center} 
\end{center}

\section{Introduction\label{section1}}

The goal of this article is to show existence of examples of discrete 
triply-periodic minimal surfaces that are modelled on smooth 
triply-periodic minimal surfaces.  For each smooth minimal surface 
model we consider, we show how one can find 
a variety of corresponding discrete minimal surfaces.  We restrict 
ourselves to discrete surfaces with a high degree of symmetry with 
respect to their density of vertices, and thus they have a highly 
discretized appearance.  The advantage of such a restriction is two-fold: 
(1) we can give explicit mathematical proofs of minimality without relying on 
numerics, and (2) we can make changes in the symmetries that would not 
be allowed in the smooth case.  In general one can also consider 
discrete minimal surfaces with finer triangulations, but considering 
only highly discretized examples is still sufficient to show existence of 
many discrete triply-periodic minimal surfaces corresponding to a single 
smooth one.  For motivational purposes, we first briefly introduce smooth 
minimal surfaces.  

\subsection{Smooth minimal surfaces}

Soap films that do not contain bounded pockets of air are surfaces that 
minimize area with respect to their boundaries.  Smooth compact minimal surfaces are 
mathematical models for soap films, because they are (by definition) 
surfaces that are critical for area with respect to all smooth variations that 
fix their boundaries.  By computing the first derivative of 
area for a smooth variation of a general surface, one finds that minimal surfaces 
are also those whose two principal curvatures at each point are equal and opposite.  
Since the mean curvature at a point on a surface is the average of the principal 
curvatures, a minimal surface is then one for which the mean curvature is 
zero at every point.  

The simplest example of a minimal surface is the flat 
plane.  Another example is the catenoid, which is a surface of revolution 
that can be parametrized by 
\begin{equation}\label{smoothcatenoid} 
\{ (\cosh x \cos y,\cosh x \sin y,x) \in \mathbb{R}^{3} \, | \, 
x \in \mathbb{R} , y \in (0,2\pi] 
\subset \mathbb{R} \} \; , \end{equation} where $\mathbb{R}$ denotes the real numbers 
and $\mathbb{R}^{3}$ denotes Euclidean $3$-space.  
(Note that by restricting to soap films not 
containing bounded pockets of air, we have ruled out surfaces like the 
sphere, which is certainly a soap film but also has nonzero mean curvature.)  

As our goal is to study discrete minimal surfaces, with smooth minimal surfaces 
playing only a suggestive role, we will not go into more detail here 
about basic properties of the smooth case.  (Some 
fine general introductions to smooth minimal surfaces are \cite{H}, 
\cite{H2}, \cite{HM}, \cite{K5}, \cite{M1}, \cite{Nit} and \cite{Oss}.)  
We will simply go directly to this definition: 

\begin{definition}
A \emph{smooth minimal surface} in ${\mathbb{R}^{3}}$ is a 
$C^\infty$ immersion $f: \mathcal{M} \to {\mathbb{R}^{3}}$ of a $2$-dimensional 
manifold $\mathcal{M}$ whose mean curvature is identically zero; or equivalently, 
the map $f$ is critical for area with respect to all 
smooth variations compactly-supported in the interior of $\mathcal M$.  
\end{definition}

\begin{figure}[tbp]
\centerline{
        \hbox{
		\psfig{figure=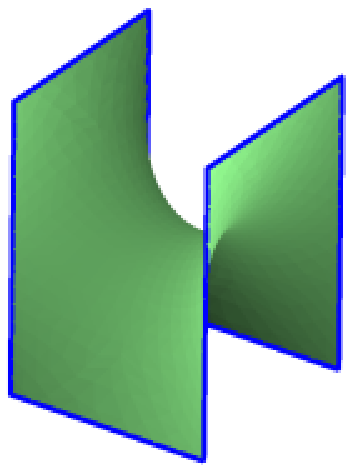,width=1.2in}
		\psfig{figure=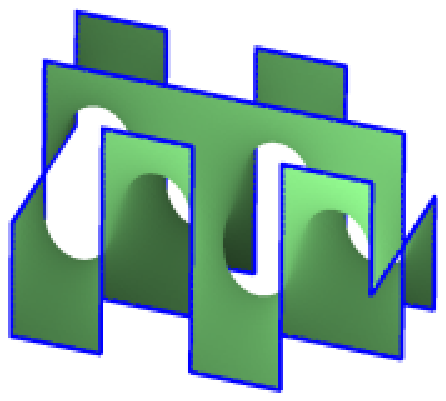,width=1.7in}
		\psfig{figure=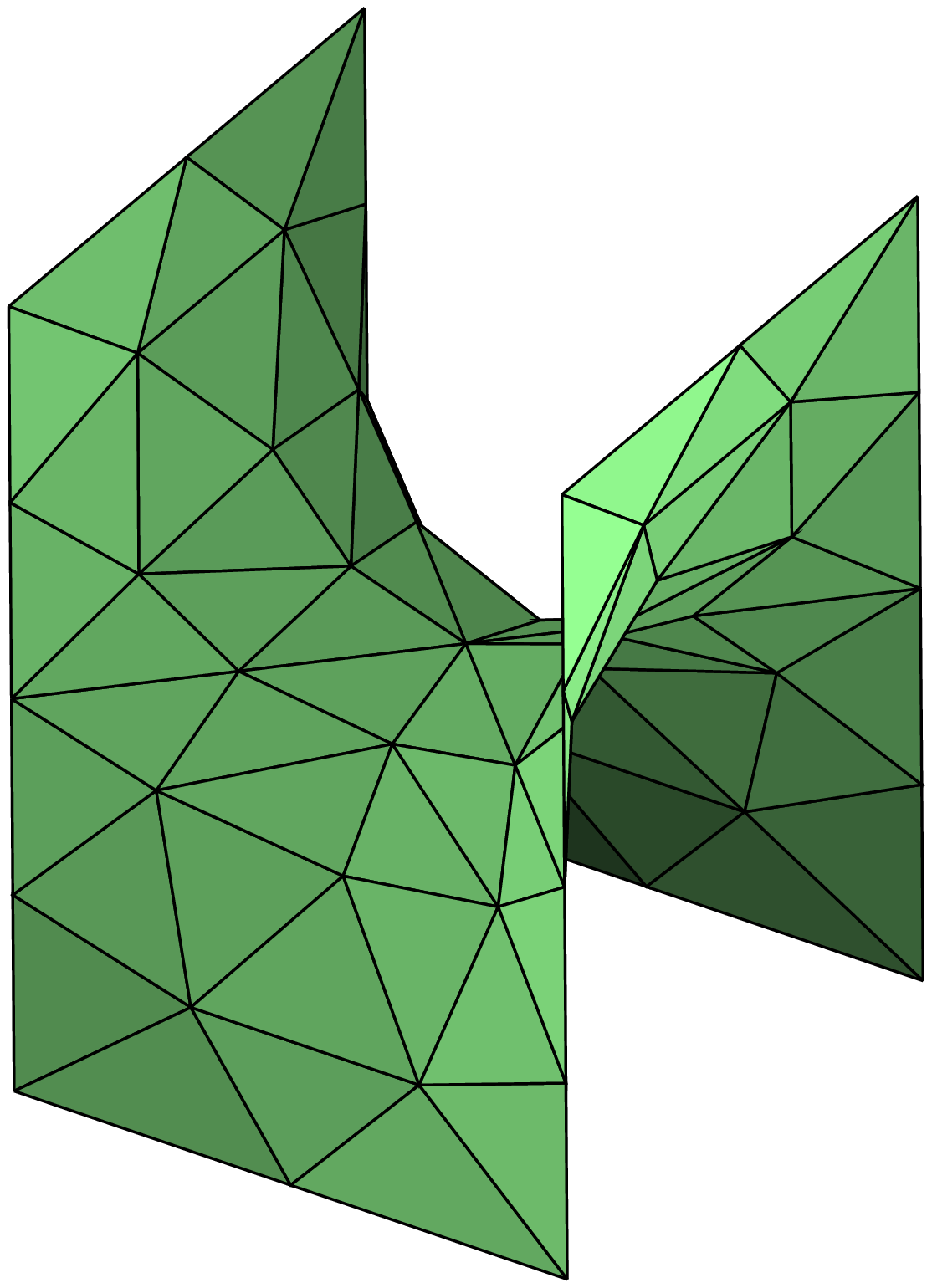,width=1.4in}
	}
}
\caption{\protect\small Smooth and discrete superman surfaces.
(All the computer graphics here were made using K. Polthier's Javaview software 
\cite{P1}.)}
\label{fig:example0}
\end{figure}

A \emph{triply-periodic} smooth 
surface is one that is periodic in three independent directions of 
${\mathbb{R}^{3}}$, i.e. $f$ and $f+\vec v_j$ have equal images for three 
independent constant vectors $\vec v_1,\vec v_2,\vec v_3$ in ${\mathbb{R}^{3}}$.  
There are a wide variety of smooth triply-periodic minimal surfaces, as can 
be seen by looking at papers of H. Karcher, K. Polthier, A. Schoen, 
W. Fischer and E. Koch \cite{FC}, \cite{K3}, \cite{K5}, \cite{KP}, \cite{Schoen}, 
for example.  We show a few examples here: the central surface 
in Figure \ref{fig:example0} was named the superman surface by W. Meeks \cite{M1}
(one special case of surfaces of this type is the Schwarz D surface); 
the second surface in the second row of Figure 
\ref{fig:example11} is a generalized Schwarz P surface (one special case of 
surfaces of this type is the original Schwarz P surface); the smooth Schwarz CLP 
surface is shown in the lower-right of Figure \ref{fig:example1}; and the 
smooth triply-periodic Fischer-Koch surface is shown in the lower row of 
Figure \ref{fig:example10}.  

\subsection{Defining discrete minimal surfaces}

Recently, finding discrete analogs of smooth objects has become 
an important theme in mathematics, appearing in a variety 
of places in analysis and geometry.  So, minimal surface theory 
being a subject in geometry that relies heavily on analysis, it is natural 
to ask how to define a discrete analog of smooth minimal surfaces.  But 
there is no single definitive answer; the definition 
one chooses would depend on which particular properties of smooth minimal surfaces 
one would wish to emulate in the discrete case.  

For example, a definition by A. Bobenko and U. Pinkall 
\cite{Bobenko/Pinkall96DiscreteIsothermic} uses 
discrete integrable systems, in analogy to smooth integrable systems properties 
of smooth minimal surfaces (or, more accurately, of smooth surfaces with 
possibly-nonzero constant mean curvature).  Their definition is good from 
the viewpoint of integrable systems, but does not yield discrete 
surfaces that are critical for area with respect to variations of their vertices.  

Here we will take a variational point of view, so we 
wish to consider area-critical discrete surfaces.  We choose this definition: 
A discrete minimal surface in ${\mathbb{R}^{3}}$ is a 
piecewise linear triangulated surface that is critical for area with respect 
to any compactly-supported boundary-fixing continuous piecewise-linear variation that 
preserves its simplicial structure, as defined in \cite{Pinkall/Polthier93Discrete} 
and \cite{PR} and Section \ref{sect.discrete} here.  

Although we will define discrete minimal surfaces as just above, 
we should remark that, even from within the variational point of 
view, this is not the unique choice of a definition.  For example, a broader 
definition by K. Polthier \cite{P2} uses "non-conforming triangulations", 
unlike the triangulations here which will all be conforming.  
Non-conforming discrete surfaces are those for which 
adjacent triangles are required only to intersect at midpoints of their 
boundary edges, not along entire edges.  Then discrete minimality can again be 
defined variationally.  This broader approach is useful for finding pairs of 
noncongruent isometric discrete minimal surfaces.  Such pairs are 
called \emph{conjugate} minimal surfaces in the case of smooth minimal surfaces.  

\subsection{Constructing discrete minimal surfaces}

Just as for smooth surfaces, a triply-periodic discrete 
surface is one that is periodic in three independent directions of 
${\mathbb{R}^{3}}$.  As noted at the beginning of this article, we will show 
that one can generally construct a variety of discrete examples modelled on a 
single smooth example.  (Hence, in this sense, there exist an 
even wider variety of discrete triply-periodic minimal surfaces than there are 
smooth ones.)  There are at least two possible ways to do this:  
\begin{itemize}
\item \underline{\bf Method 1} vary the choice of simplicial structure of some compact 
portion of the surface, and 
\item \underline{\bf Method 2} vary the choice of rigid motions of 
${\mathbb{R}^{3}}$ that create the complete surface from some compact portion.
\end{itemize}
To explain these two methods in more detail, imagine a compact portion $M$ of 
a smooth triply-periodic minimal surface with piece-wise smooth boundary 
$\partial M$ consisting of smooth curves $\gamma_1$,...,$\gamma_n$.  Suppose 
that each $\gamma_j$ is 
either a straight line segment or a curve in a principal curvature direction of $M$ 
that also lies in a plane of ${\mathbb{R}^{3}}$.  (In the latter case we call 
$\gamma_j$ a \emph{planar geodesic}, since it is necessarily a geodesic 
of $M$.)  A larger minimal 
surface $M^\prime$ is constructed from $M$ by including images of $M$ under 
$180^\circ$ rotations 
about the lines containing linear $\gamma_j$ and under reflections 
through the planes containing planar geodesic $\gamma_j$.  (The fact that the 
larger portion $M^\prime$ is still a smooth minimal surface can be shown using 
complex analysis, see \cite{K3}, \cite{K5}, \cite{Nit}, \cite{Oss}, 
for example.)  The larger portion 
$M^\prime$ again has a piece-wise smooth boundary $\partial M^\prime$ 
consisting of line segments and planar geodesics, so this procedure can be 
repeated again on $M^\prime$.  Repeating this 
procedure on ever-bigger pieces of the surface a countably infinite 
number of times, one builds the entire complete surface.  
$M$ is often called a \emph{fundamental domain} of the complete surface.  

For example, the minimal surface on the left-hand side of Figure 
\ref{fig:example0} is a fundamental domain $M$ of a complete triply-periodic 
smooth minimal surface.  The boundary $\partial M$ contains the eight vertices 
\[ p_1=(0,0,0) , \;\; p_2=(x,0,0) , \;\; p_3=(x,0,z) , \;\; p_4=(x,y,z) , \]\[ 
p_5=(x,y,0) , \;\; p_6=(0,y,0) , \;\; p_7=(0,y,z) , \;\; p_8=(0,0,z) , \] 
for some given positive reals $x,y,z>0$.  
Then $\partial M$ is a polygonal loop consisting of eight line segments, 
from $p_j$ to $p_{j+1}$ for $j=1,2,...,7$, and 
finally from $p_8$ to $p_1$.  One can construct the entire complete 
surface, as described above, using only $180^\circ$ rotations about boundary 
line segments.  A larger piece of this complete surface can be seen in 
the center of Figure \ref{fig:example0}.  For general values of 
$x,y,z$, the resulting complete triply-periodic surface is a 
superman surface.  When $x=y$ this surface represents 
Schwarz' solution of Gergonne's problem (see \cite{K3} and \cite{KP} for more 
on this), and when $x=y=\sqrt{2} \cdot z$ this surface is the Schwarz D surface.  

\begin{figure}[tbp]
\centerline{
        \hbox{
		\psfig{figure=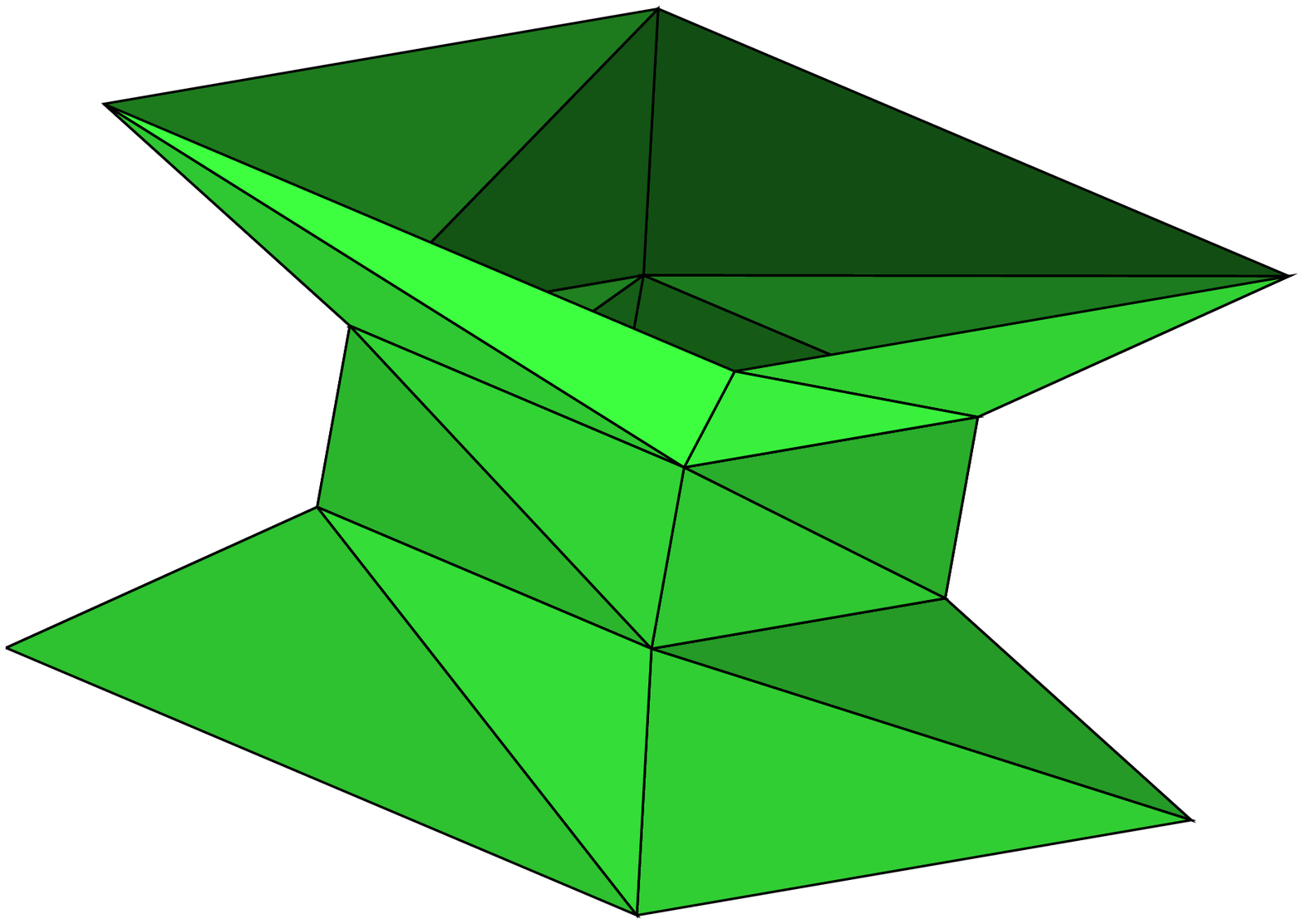,width=1.0in}
		\psfig{figure=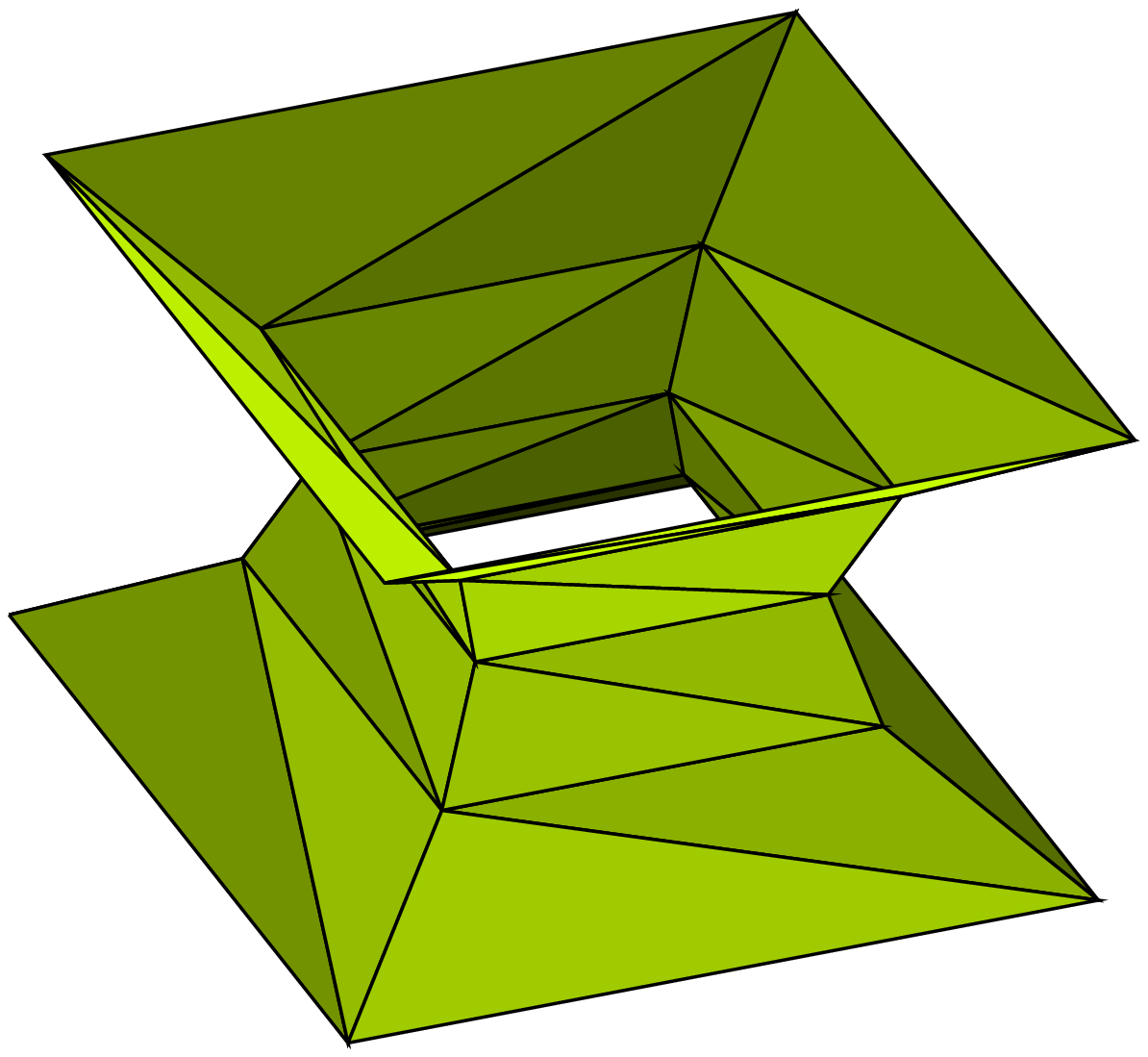,width=1.0in}
		\psfig{figure=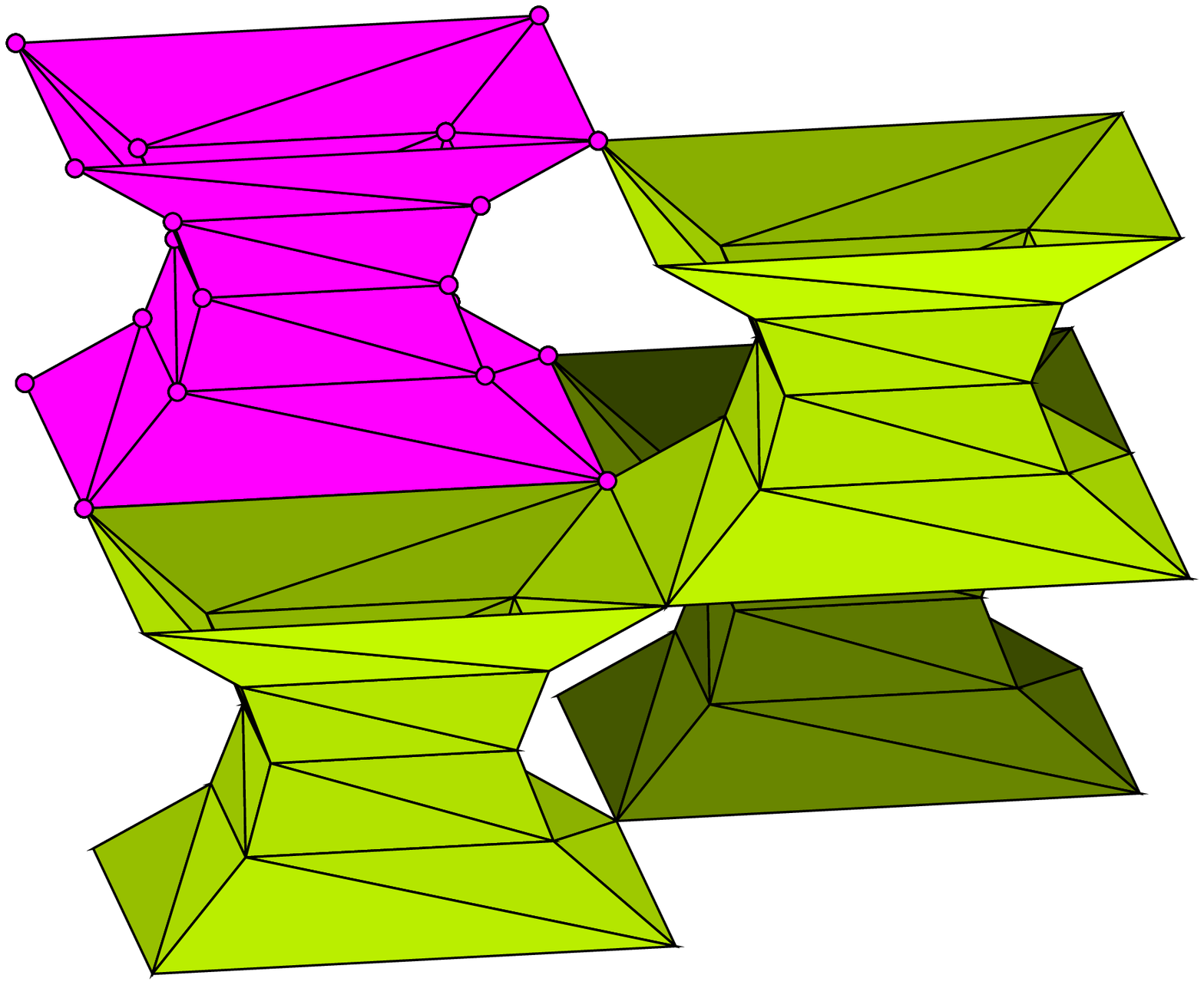,width=1.4in}
		\psfig{figure=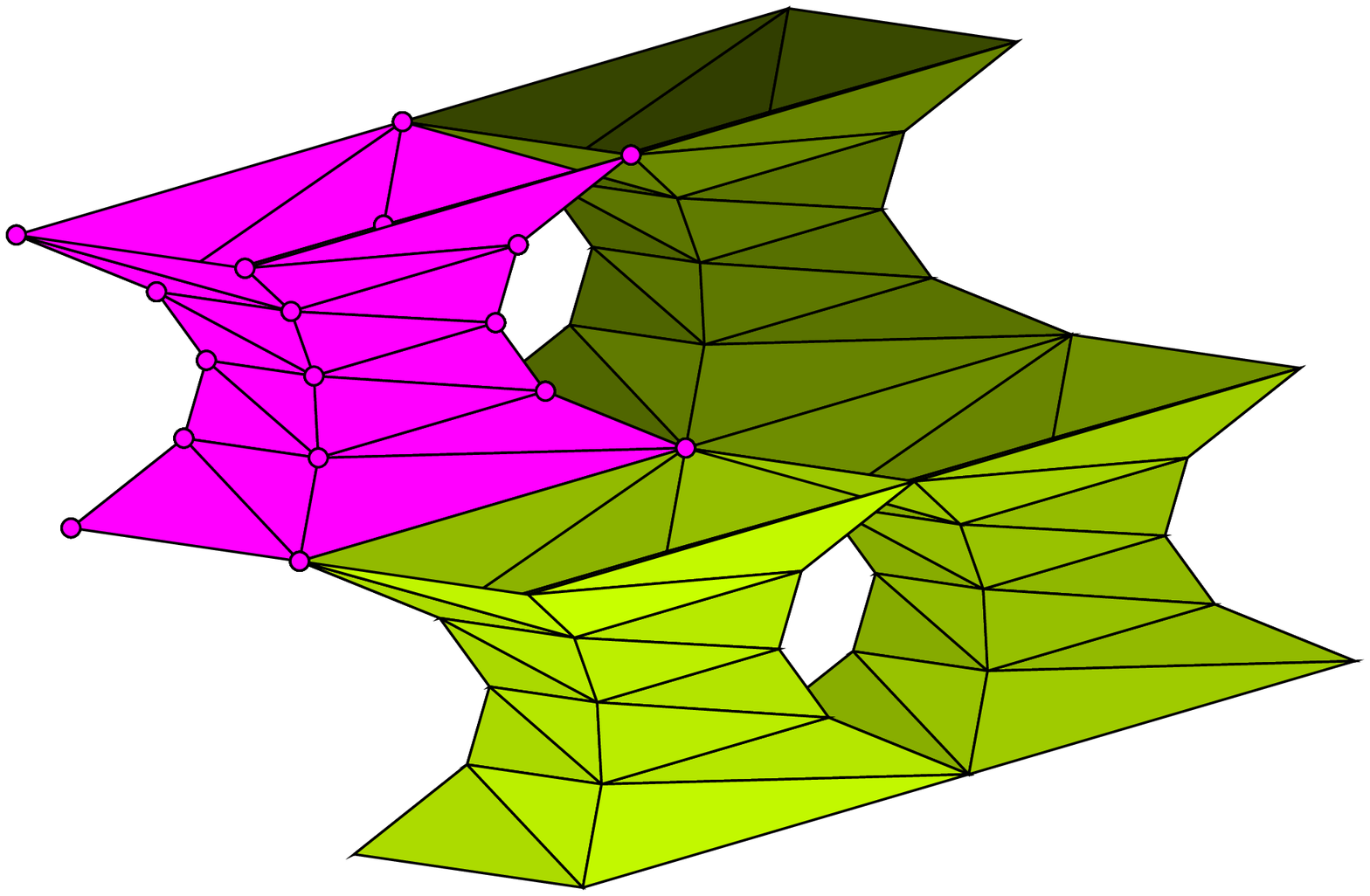,width=1.4in}
	}
}
\centerline{
        \hbox{
		\psfig{figure=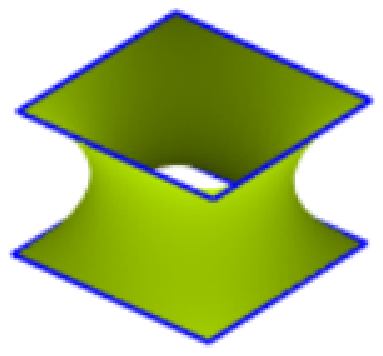,width=1.0in}
		\psfig{figure=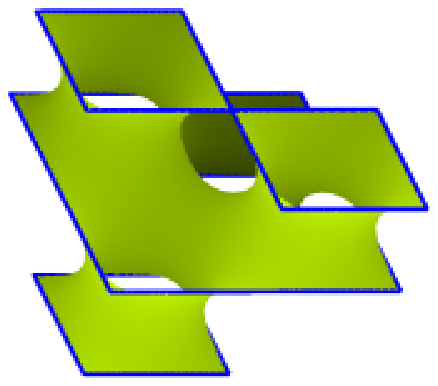,width=1.6in}
		\psfig{figure=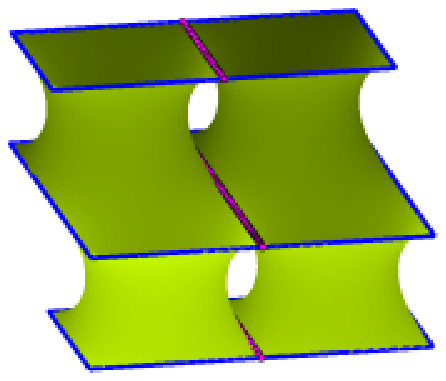,width=1.6in}
		\psfig{figure=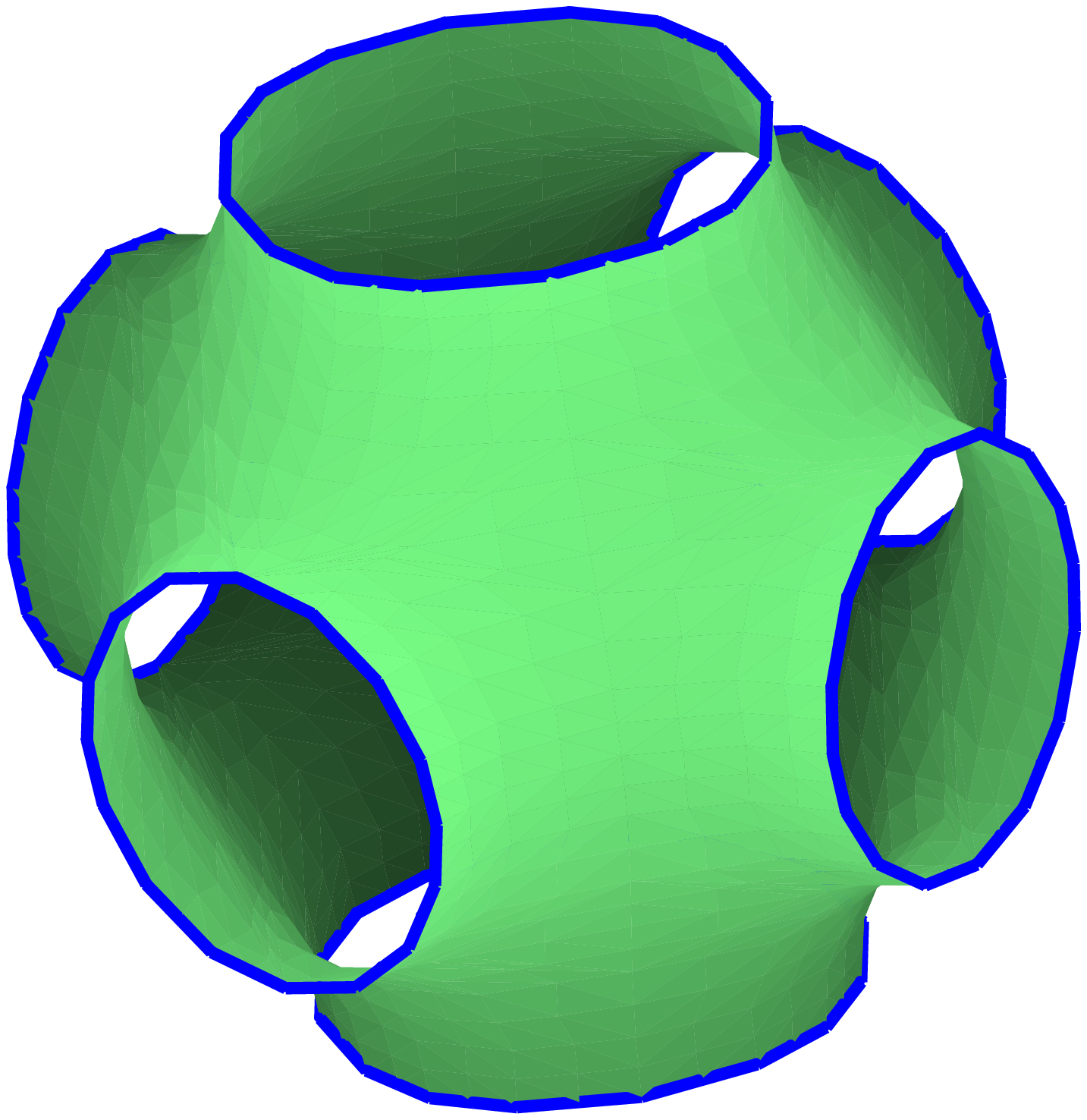,width=1.0in}
	}
}
\caption{\protect\small Discrete and smooth Schwarz P surfaces.}
\label{fig:example11}
\end{figure}

To construct triply-periodic discrete minimal surfaces modelled on smooth superman 
surfaces, we first find a discrete version of the smooth fundamental domain $M$.  One 
such example, found numerically using JavaView software \cite{P1}, 
is on the right-hand side of Figure \ref{fig:example0}.  However, there are many 
ways the simplicial structure of the discrete version can be chosen, and a 
number of them are shown in Figure \ref{fig:example8}.  
(The examples in Figure \ref{fig:example8} 
all have coarser simplicial structures, and hence showing their minimality is 
manageable by direct computation without using a computer, as we will 
see in Section \ref{section3}.)  The corresponding 
complete triply-periodic discrete minimal surfaces are then constructed in 
the same way as in the smooth case, by $180^\circ$ 
rotations about boundary line segments.  This variety of simplicial structures all 
based on one smooth superman surface is an example of using {\bf Method 1}.  

To demonstrate {\bf Method 2}, consider the left-hand 
minimal surface $M$ in the second row of Figure \ref{fig:example11}.  Its 
boundary $\partial M$ is two squares in parallel planes, where 
one square projects to the other by projection 
orthogonal to the planes.  As in the previous example, $180^\circ$ rotations about 
boundary line segments produces a complete 
triply-periodic smooth minimal surface, and a larger piece of this surface 
is shown just to the right of $M$ in Figure \ref{fig:example11}.  
This example is a Schwarz P surface.  
The original Schwarz P surface is the special case that the distance between 
the two parallel planes containing $\partial M$ is 
$\sqrt{2}$ times shorter than the length of each edge in 
$\partial M$, and is conjugate to the Schwarz D surface.  
(For later use, we mention that the lower-right surface $\hat M$ in Figure 
\ref{fig:example11} is also a portion of a Schwarz P surface, now 
bounded by six planar geodesics, each lieing within one face of a 
rectanguloid with square base.  The entire complete surface can be built 
from $\hat M$ solely by applying a set of reflections in equally-spaced planes 
parallel to the faces of the rectanguloid.  The top one-fourth of $\hat M$ 
equals the bottom half of $M$.  So one could choose either $M$ or $\hat M$ as 
the fundamental piece for constructing a complete Schwarz P surface.)  

As with the superman surface, to create discrete analogs 
of the Schwarz P surface, one can apply {\bf Method 1} and choose amongst many 
different simplicial structures for this surface.  
Two such possibilities are the first two discrete minimal surfaces in 
the first row of Figure \ref{fig:example11}.  The first one has $4$ 
squares in parallel planes amongst its edge set, and the second one has 
$5$ squares in parallel planes amongst its edge set.  In fact, one can 
make examples with any number$\geq 3$ of such squares in its edge set, as 
we will see in Section \ref{section4}.  In the third figure 
in the first row of Figure \ref{fig:example11}, 
a larger portion of a resulting complete discrete triply-periodic minimal 
surface is shown.  

But let us return to a demonstration of {\bf Method 2} on the Schwarz P 
surface.  Let $P$ be a plane perpendicular 
to the planes containing $\partial M$ so that $P$ also contains two disjoint 
boundary edges in $\partial M$.  Consider the surface we would 
have if we first included a reflected image of $M$ across $P$, and then 
created a complete surface by $180^\circ$ rotations about all 
resulting boundary line segments.  A part of this surface is shown in the 
third figure of the second row of Figure \ref{fig:example11} (this part is 
bounded by three boundary components, one square and two rectangles).  In this 
part in Figure \ref{fig:example11}, there are three parallel dotted lines, 
along which the surface is not even a $C^1$ immersion, hence the mean curvature 
is not defined there, and we cannot call this a minimal surface.  So this 
construction is forbidden in the smooth case.  
However, for the discrete analogs shown in Figure \ref{fig:example11}, such 
a construction actually does produce a discrete 
complete triply-periodic minimal surface.  This construction works in the 
discrete case because the notion of "minimality" is defined \emph{only} at the 
vertices of the discrete surface, as we will see in Section \ref{sect.discrete}.  
Choosing the simplicial structure so that there are no vertices in 
the interiors of the edges comprising $\partial M$ is what makes this construction 
possible.  

Thus we have varied 
the rigid motions of ${\mathbb{R}^{3}}$ that were used to create 
the complete smooth surface from $M$, in a way that cannot be allowed 
for smooth surfaces, to make a different kind of discrete minimal surface 
modelled on $M$.  Part of this surface is shown in the upper-right of 
Figure \ref{fig:example11} (this part also has three boundary components, one 
square and two rectangles).  One can easily imagine including more reflections (not 
just across a single plane $P$), to make infinitely many different 
kinds of discrete minimal surfaces modelled on $M$.  This is {\bf Method 2}.  

In Section \ref{mainresults}, we state our main result.  
In Sections \ref{section3} and \ref{section4}, as applications of {\bf Method 1}, 
we will see a variety of different simplicial 
structures making discrete analogs of the smooth superman and 
Schwarz P surfaces.  Section \ref{section5} contains examples 
modelled on other smooth minimal surfaces.  In the last example 
of Section \ref{section4} and the first example of Section \ref{section5} there are 
further applications of {\bf Method 2}.  

\section{Discrete Minimal Surfaces\label{sect.discrete}}

In the smooth case, a compact minimal surface is area-critical for
any variation that fixes the boundary. We wish to define discrete minimal 
surfaces so that they have the same variational property for the same types of 
variations.  We begin by defining discrete surfaces and their variations, and we 
first give an informal definition: 

\begin{definition} (Informal Definition) 
A \emph{discrete surface} in ${\mathbb{R}^{3}}$ is a 
$C^0$ mapping $f: \mathcal{M} \to {\mathbb{R}^{3}}$ of a $2$-dimensional 
manifold $\mathcal{M}$ so that each face of some triangulation of 
$\mathcal{M}$ is mapped to a triangle in ${\mathbb{R}^{3}}$.  The surface 
$f(\mathcal{M})$ is \emph{embedded} if $f$ is an injection.  
\end{definition}

We define \emph{embedded} in the discrete case without any conditions about 
nondegeneracy of $f$ (nondegeneracy is meaningless here, as $f$ is only 
$C^0$).  However, we still use this word, to maintain the 
analogy to embeddedness of smooth surfaces.  

We stated the above definition for the intuition it provides, but 
we will require a more involved definition, which we now give: 

\begin{definition} (Formal Definition) 
A \emph{discrete surface} in ${\mathbb{R}^{3}}$ is a triangular mesh 
which has the topology of an abstract $2$-dimensional 
locally-finite simplicial surface $K$ combined with a geometric $C^0$ 
realization $\mathcal T$ in ${\mathbb{R}^{3}}$ that is 
piecewise-linear on each simplex.  (Because $K$ is a simplicial "surface", each 
$1$-dimensional simplex in $K$ lies in the boundary of exactly one or 
two $2$-dimensional simplices of $K$.)  
The geometric realization $\mathcal T$ is determined by a set of vertices $\mathcal{%
V} = \{p_1,p_2,...\} \subset {\mathbb{R}^{3}}$ corresponding to the $0$-dimensional 
simplices of $K$, and $\mathcal V$ could 
be either finite or countably infinite.  The simplicial surface $K$
represents the connectivity of $\mathcal T$. The $0$, $1$, and $2$
dimensional simplices of $K$ represent the vertices, edges, and triangles of
$\mathcal T$.  

\begin{figure}
\begin{center}
\unitlength=0.5pt
\begin{picture}(200.00,200.00)(0.00,0.00)
\put(70,130){\circle*{3}}
\put(110,40){\circle*{3}}
\put(190,110){\circle*{3}}
\put(180,170){\circle*{3}}
\put(120,190){\circle*{3}}
\put(130,130){\circle*{3}}
\bezier200(70,130)(90,85)(110,40)
\bezier200(190,110)(150,75)(110,40)
\bezier100(190,110)(185,140)(180,170)
\bezier100(120,190)(150,180)(180,170)
\bezier100(120,190)(95,160)(70,130)
\put(85,94){\vector(3,2){43}}
\put(150,75){\vector(-1,1){42}}
\put(185,140){\vector(-4,-1){30}}
\put(150,180){\vector(-2,-3){22}}
\put(95,160){\vector(1,-1){23}}
\bezier100(70,130)(95,130)(130,130)
\bezier200(110,40)(120,85)(130,130)
\bezier100(190,110)(160,120)(130,130)
\bezier100(180,170)(155,150)(130,130)
\bezier100(120,190)(125,160)(130,130)
\put(131,116){$p$}
\put(195,105){$q$}
\put(185,175){$r$}
\put(187,140){$J(r-q)$}
\put(108,150){$\mathcal{T}$}
\end{picture}
\end{center}
\vspace{-0.2in}
\caption{\protect\small
At each vertex $p$ the gradient of discrete area is the sum of the 
$90^\circ$-rotated edge
vectors $J(r-q)$, as in Equation \eqref{eq.gradArea}.}\label{basicfigure}
\end{figure}
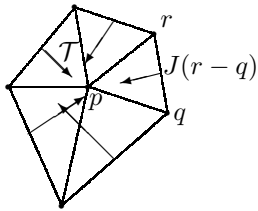

Let $T=(p,q,r)$ denote an oriented triangle of $\mathcal{T}$ with vertices $%
p,q,r \in \mathcal{V}$. Let $\overline{pq}$ denote an edge of $T$ with
endpoints $p,q \in \mathcal{V}$.

For $p\in \mathcal{V}$, let ${\func{star}}(p)$ denote the triangles of $%
\mathcal{T}$ that contain $p$ as a vertex. 

The \emph{boundary} $\partial {\mathcal T}$ of 
$\mathcal T$ is the union of those edges 
bounding only a single triangle of $\mathcal T$.  
The \emph{interior vertices} (respectively, \emph{boundary vertices}) of 
$\mathcal T$ are those that are not contained (respectively, are 
contained) in $\partial {\mathcal T}$.  

We say that $\mathcal T$ is \emph{complete} if $\partial \mathcal{T}$ is empty and 
if $\mathcal T$ is complete with respect to the distance function induced by its 
realization in ${\mathbb{R}^{3}}$.  
\end{definition}

\begin{definition}\label{defnofvariationsurfaces}
Let $\mathcal{V}=\{p_{1},p_{2}, ... \}$ be the set of vertices of a
discrete surface $\mathcal{T}$. A \emph{variation} $\mathcal{T}(t)$ of $%
\mathcal{T}$ is defined as a $C^{\infty}$ variation of the vertices $p_{i}$ 
\[ p_{i}(t):[0,\epsilon)\rightarrow {\mathbb{R}^{3}} \mbox{  so that  }
p_{i}(0)=p_{i}\;\forall i=1,...,m. 
\]
The straightness of the edges and the 
flatness of the triangles are preserved as the vertices $p_i(t)$ move 
with respect to $t$.  

When $\mathcal T$ is compact, we say that ${\mathcal T}(t)$ 
\emph{fixes the boundary} $\partial {\mathcal T}$ if $p_i(t)$ is constant 
in $t$ for all $p_i \in \partial \mathcal{T}$.  
When $\mathcal T$ is complete, we say that ${\mathcal T}(t)$ 
is \emph{compactly supported} if $p_i(t)$ is constant in $t$ for all 
but a finite number of vertices $p_i$.  
\end{definition}

The area of a discrete surface is 
\[
\func{area} \mathcal{T}:=\sum_{T\in \mathcal{T}}\func{area}T \; , 
\]
where $\func{area}T$ denotes the Euclidean 
area of the triangle $T$ as a subset of ${{\mathbb{R}^{3}}}$.

\begin{lemma}\label{firstlemma}
Let $\mathcal{T}(t)$ be a variation of a discrete surface $\mathcal{T}$. At
each vertex $p$ of $\mathcal{T}$, the gradient of area is 
\begin{equation}
\nabla_{p}\func{area}\mathcal{T}={\frac{1}{2}}\sum_{{\ T=(p,q,r) \in \func{%
star}(p)}}J(r-q) \; ,  \label{eq.gradArea} 
\end{equation}
where $J$ is $90^\circ$ rotation in the plane of each
oriented triangle $T$. The first derivative of the area is then
given by the chain rule 
\begin{equation}\label{firstderivativeofarea} \left. 
\frac{d}{dt}\func{area}\mathcal{T}(t) \right|_{t=0}=\sum_{p\in \mathcal{V}} 
\left\langle \left. \frac{d(p(t))}{dt} 
\right|_{t=0},\nabla_{p}\func{area}\mathcal{T} \right\rangle \; . 
\end{equation}
\end{lemma}

\begin{proof}
Let $p_i(t)$ be the corresponding variation of each vertex in the vertex 
set ${\mathcal V}(t)$ of the variation ${\mathcal T}(t)$.  Then 
\[ \func{area}(\mathcal{T}(t)) = \frac{1}{6} \sum_{p(t) \in {\mathcal V}(t)} 
\left( \sum_{(p(t),q(t),r(t)) \in \func{star}(p(t))} 
|| (r(t)-p(t)) \times (q(t)-p(t)) || \right) \; , \] and a computation 
implies \[ \frac{d}{dt} \func{area}(\mathcal{T}(t)) = \frac{1}{2} 
\sum_{p(t) \in {\mathcal V}(t)} \left\langle \frac{d(p(t))}{dt} , 
\sum_{(p(t),q(t),r(t)) \in \func{star}(p(t))} 
||r(t)-q(t)|| \eta(t) \right\rangle \; , \] 
where $\eta(t)$ is the unit conormal in the plane of the triangle 
$(p(t),q(t),r(t))$ along the edge $r(t)-q(t)$, oriented in the 
same direction as $J(r(t)-q(t))$.  Restricting to $t=0$ proves the lemma.  
\end{proof}

As defined in Section \ref{section1}, 
a smooth immersion $f: \mathcal{M} \to {\mathbb{R}^{3}}$ of a $2$-dimensional 
complete manifold $\mathcal M$ without boundary is minimal 
if $f$ is area-critical for all 
compactly-supported smooth variations.  In the case that $\mathcal M$ 
is compact with boundary, then $f$ is minimal if it is 
area-critical for all smooth variations preserving 
$f(\partial {\mathcal M})$.  

We wish to define discrete minimal surfaces $\mathcal T$ so that they have 
the analogous properties, for variations as in Definition 
\ref{defnofvariationsurfaces}.  
So when $\mathcal T$ is compact, we consider variations $\mathcal{T}(t)$ of 
$\mathcal{T}$ that fix $\partial \mathcal{T}$; 
and when $\mathcal T$ is complete, we consider variations 
$\mathcal{T}(t)$ of $\mathcal{T}$ that are compactly supported.  
By Lemma \ref{firstlemma}, the condition that makes $\mathcal T$ 
area-critical for any 
variation of these types is expressed in the following definition.  

\begin{definition}\label{cmcdefinition}
A discrete surface is \emph{minimal} if \begin{equation}\label{eqn:gradientiszero} 
\nabla_{p}\func{area} {\mathcal T}=0 \end{equation} for all interior vertices $p$. 
\end{definition}

\begin{remark}\label{flatpartsofminimalsurfaces}
If $\mathcal{T}$ is a discrete minimal surface that contains a 
discrete subsurface $\mathcal{T}^\prime$ lieing in a plane $P$, it follows 
from Equations \eqref{eq.gradArea} and \eqref{eqn:gradientiszero} that the 
discrete minimality of 
$\mathcal{T}$ is independent of the choice of triangulation of the trace of 
$\mathcal{T}^\prime$ within $P$.  Thus whenever such a planar part 
$\mathcal{T}^\prime$ 
occurs in the following examples, we will be free to triangulate 
$\mathcal{T}^\prime$ any way we please, within its trace in $P$.  
\end{remark}

\section{Results\label{mainresults}}

For the purpose of stating our main theorem, we give the following two 
definitions: 

\begin{definition}
A discrete triply-periodic minimal surface $\mathcal T$ has 
{\em common topology and symmetry} as a smooth triply-periodic minimal 
immersion $f: {\mathcal M} \to {\mathbb{R}^{3}}$ if there exists a homeomorphism 
\[ \phi : f(\mathcal{M}) \to \mathcal{T} \] such that the following statement 
holds: $R_s : 
{\mathbb{R}^{3}} \to {\mathbb{R}^{3}}$ is a rigid motion preserving $f(\mathcal{M})$ 
if and only if there exists a rigid motion $R_d : 
{\mathbb{R}^{3}} \to {\mathbb{R}^{3}}$ preserving $\mathcal{T}$ so that 
\[ R_d \circ \phi = \phi \circ R_s |_{f(\mathcal{M})} \; , \]  
and furthermore $R_s$ is a reflection (resp. translation, rotation, screw motion) 
if and only if $R_d$ is a reflection (resp. translation, rotation, screw motion).  
\end{definition}

\begin{definition}
We say that a subsurface $\mathcal{T}^\prime$ of a complete 
discrete triply-periodic minimal surface $\mathcal T$ is a 
{\em fundamental domain} if $\mathcal{T}^\prime$ can be extended to all of $\mathcal{T}$ 
by a discrete group of rigid motions $\{R_{d,\alpha} \}_{\alpha \in \Lambda}$ generated by 
\begin{enumerate}
\item reflections across planes containing boundary edges and 
\item $180^\circ$ degree rotations about boundary edges
\end{enumerate}
so that each $R_{d,\alpha}$ is a symmetry of the full surface $\mathcal{T}$.  
\end{definition}

\begin{remark}
In the above definition of a fundamental domain, we do not allow 
rigid motions that do not fix any edges of $\mathcal{T}$ (thus any 
fundamental domain of the example in Subsection \ref{section3-1} must contain at least $6$ 
triangles, even in the most symmetric case $x=1$).  Also, we do not allow rigid 
motions that are not symmetries of the full surface $\mathcal{T}$ (thus any 
fundamental domain of the example in Subsection \ref{section4-2} must 
contain at least $32$ triangles).  
\end{remark}

We now state our results about embedded triply-periodic discrete 
minimal surfaces, which involve comparisons to the following smooth minimal surfaces: 
the superman surfaces (Figure \ref{fig:example0}), 
the Schwarz P surfaces (Figure \ref{fig:example11}), 
the Schwarz H surfaces, the Schwarz CLP surfaces (Figure \ref{fig:example1}), 
A. Schoen's I-Wp and F-Rd and H-T surfaces, and the triply-periodic Fischer-Koch 
surfaces (Figure \ref{fig:example10}).  
More complete information about these smooth surfaces can be found in 
\cite{FC}, \cite{H}, \cite{K3}, \cite{K5}, 
\cite{KP}, \cite{M1}, \cite{M2}, \cite{MRR}, 
\cite{Nit}, \cite{MartyRoss} and \cite{Schoen}.  

\begin{theorem}
The following discrete embedded triply-periodic minimal surfaces exist: 
\begin{enumerate}
\item those with common topology and symmetry as smooth superman surfaces whose 
fundamental domains contain $4$, $5$, $6$ or $8$ triangles; 
\item those with common topology and symmetry as smooth Schwarz P surfaces whose 
fundamental domains contain $1$, $2$, $6$ or $32$ triangles, and also a different class of 
discrete surfaces with common topology and symmetry as smooth Schwarz P surfaces whose 
fundamental domains contain $2n$ triangles for any positive integer $n$; 
\item those with common topology and symmetry as smooth Schwarz H surfaces whose 
fundamental domains contain $2n$ triangles for any positive integer $n$; 
\item those with common topology and symmetry as smooth Schwarz CLP surfaces whose 
fundamental domains contain $6$ triangles; 
\item one with common topology and symmetry as A. Schoen's smooth I-Wp surface whose 
fundamental domain contains $5$ triangles; 
\item one with common topology and symmetry as A. Schoen's smooth F-Rd surface whose 
fundamental domain contains $3$ triangles; 
\item those with common topology and symmetry as A. Schoen's smooth H-T surfaces whose 
fundamental domains contain $6$ triangles; 
\item those with common topology and symmetry as the smooth triply-periodic surfaces of 
Fischer-Koch whose fundamental domains contain $8$ triangles.  
\end{enumerate}
\end{theorem}

To prove this theorem, we need only collect the examples proven 
to exist in the remainder of this paper, as follows: 

\begin{proof}
Embedded discrete superman surfaces whose 
fundamental domains contain $4$ (resp. $5$, $6$, $8$) triangles are given in 
Subsection \ref{section3-3} (resp. \ref{section3-4}, 
\ref{section3-2}, \ref{section3-1}).  A second type of embedded 
discrete superman surfaces whose 
fundamental domains contain $6$ triangles are given in 
the second to the last paragraph of Subsection \ref{section5-1} by Method 2.  

Embedded discrete Schwarz P surfaces whose 
fundamental domains contain $1$ (resp. $2$, $6$, $32$) triangles are given in the 
first example of Subsection \ref{section4-1} (resp. the second example of Subsection 
\ref{section4-1}, Subsection \ref{section4-2}, the 
last paragraph of Subsection \ref{section5-1} by Method 2).  
Other classes of embedded discrete Schwarz P surfaces are given by 
choosing $k=4$ and $z_0=0$ and $j_0=n$ in Subsection \ref{section4-3} 
for any positive integer $n$, and here we 
can choose the fundamental domains to be the $2n$ 
triangles between two adjacent meridans and below the plane 
$\{ (x,y,0) \, | \,  x,y \in \mathbb{R} \}$.  

Embedded discrete Schwarz H surfaces are given by 
choosing $k=3$ and $z_0=0$ and $j_0=n$ in Subsection \ref{section4-3} 
for any positive integer $n$, and here again we 
can choose the fundamental domains to be the $2n$ 
triangles between two adjacent meridans and below the plane 
$\{ (x,y,0) \, | \,  x,y \in \mathbb{R} \}$.  

Embedded discrete Schwarz CLP surfaces whose 
fundamental domains contain $6$ triangles are given in Subsection \ref{section5-1}.  

Embedded discrete I-Wp and F-Rd surfaces whose 
fundamental domains contain $5$ and $3$ triangles, respectively, 
are given in Subsection \ref{section5-2}.  

Embedded discrete H-T surfaces whose 
fundamental domains contain $6$ triangles are given in Subsection \ref{section5-3}.  

Embedded discrete triply-periodic 
Fischer-Koch surfaces whose fundamental domains contain $8$ triangles 
are given in Subsection \ref{section5-4}.  
\end{proof}

\section{Discrete versions of the superman surface\label{section3}}

In Sections \ref{section3}, \ref{section4} and \ref{section5}, we construct 
discrete triply-periodic 
minimal surfaces.  All of the surfaces we construct are embedded.  

To construct examples, we always start with a compact discrete 
fundamental piece $\mathcal T$, with given simplicial structure and 
boundary constraints.  
The complete triply-periodic discrete surface is then formed by including 
images of $\mathcal T$ under a discrete group of rigid motions of 
${\mathbb{R}^{3}}$.  This group of rigid motions is generated by a finite 
number of $180^\circ$ rotations about lines and/or reflections across planes, 
and for each edge $\overline{p q}$ in $\partial 
{\mathcal T}$ this group contains either 
\begin{itemize}
\item the $180^\circ$ rotation about the line containing $\overline{p q}$, or 
\item a reflection across a plane containing $\overline{p q}$.  
\end{itemize}
To ensure that the resulting complete discrete triply-periodic surface is minimal, 
Section \ref{sect.discrete} gave us the following two approaches: 
\begin{enumerate}
\item Use symmetries of $\mathcal T$ and of the resulting complete discrete surface 
to show that Equation \eqref{eqn:gradientiszero} holds at the vertices.  
\item Locate the vertices of $\mathcal T$ so that $\mathcal T$ is area-critical 
with respect to its boundary constraints.  
\end{enumerate}
In the following examples, 
either approach produces the same conditions for minimality.  

As noted in the introduction, we wish to show examples here for which 
explicit mathematical proofs of minimality are still manageable (without 
the aid of a computer).  So we are limited to examples with a high degree of 
symmetry with respect to their density of vertices, and thus with a highly 
discretized appearance.  Discrete minimal surfaces that appear more 
like approximations of smooth minimal surfaces usually can only be found 
numerically.  Numerical examples, with finer simplicial structures, of discrete 
versions of the superman, Schwarz P, F-Rd, I-Wp and H-T surfaces are shown 
in \cite{P2}.  

\subsection{First example\label{section3-1}} 

\begin{figure}[tbp]
\centerline{
        \hbox{
		\psfig{figure=exa04a.ps,width=1.7in}
		\psfig{figure=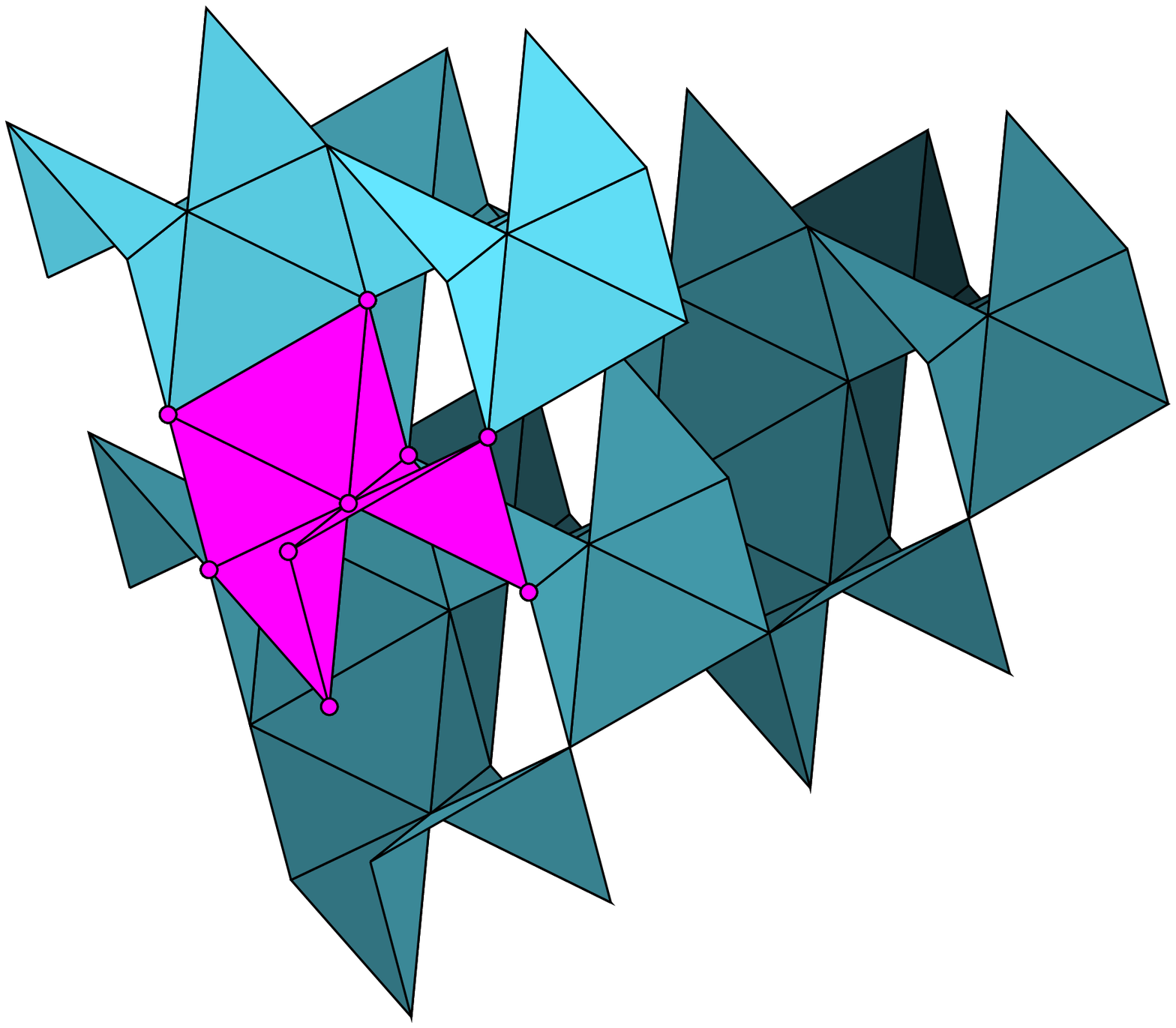,width=1.7in}
		\psfig{figure=exa04b.ps,width=1.7in}
	}
}
\centerline{
        \hbox{
		\psfig{figure=exa05a.ps,width=1.0in}
		\psfig{figure=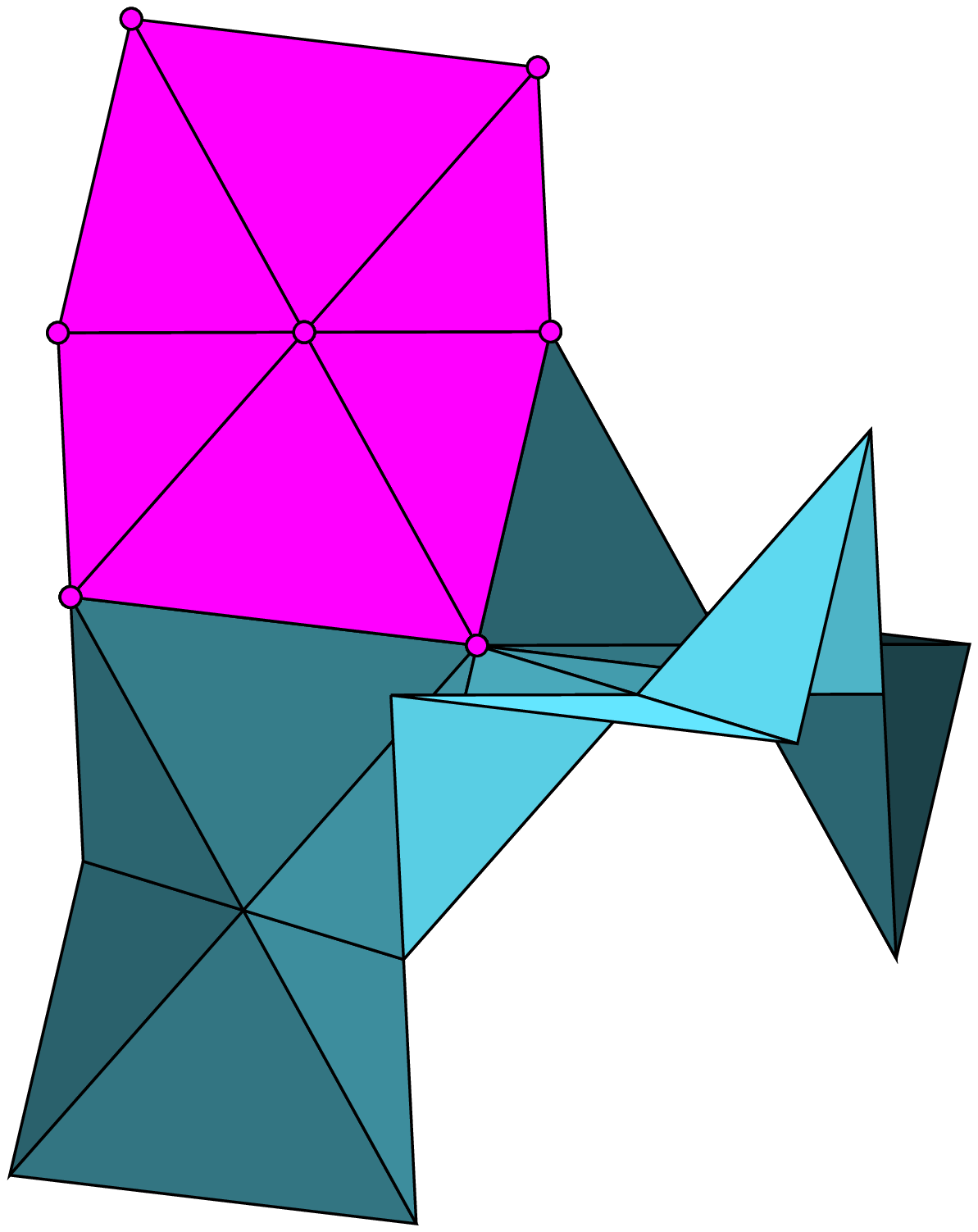,width=1.0in}
		\psfig{figure=exa05b.ps,width=1.0in}
		\psfig{figure=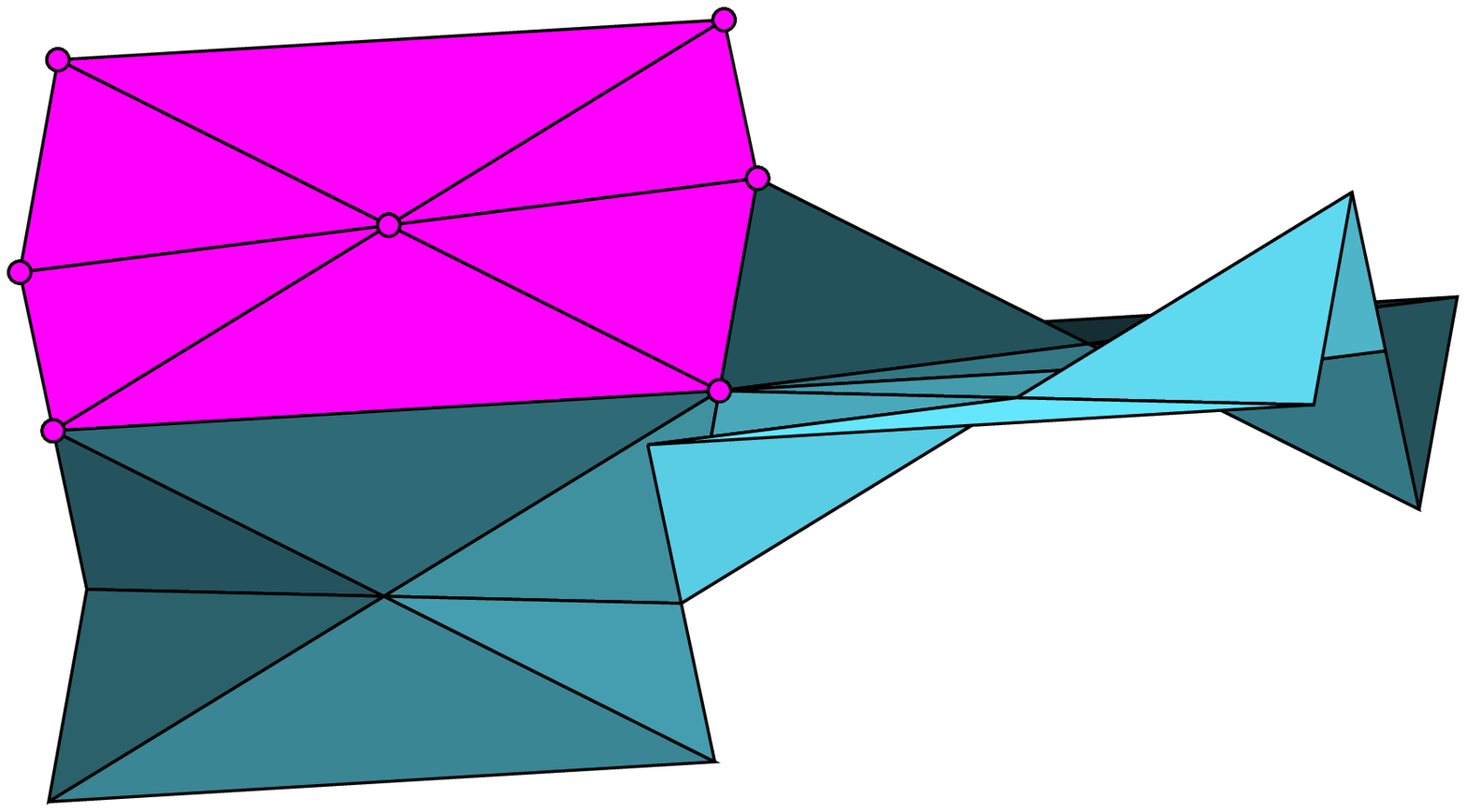,width=1.7in}
	}
}
\centerline{
        \hbox{
		\psfig{figure=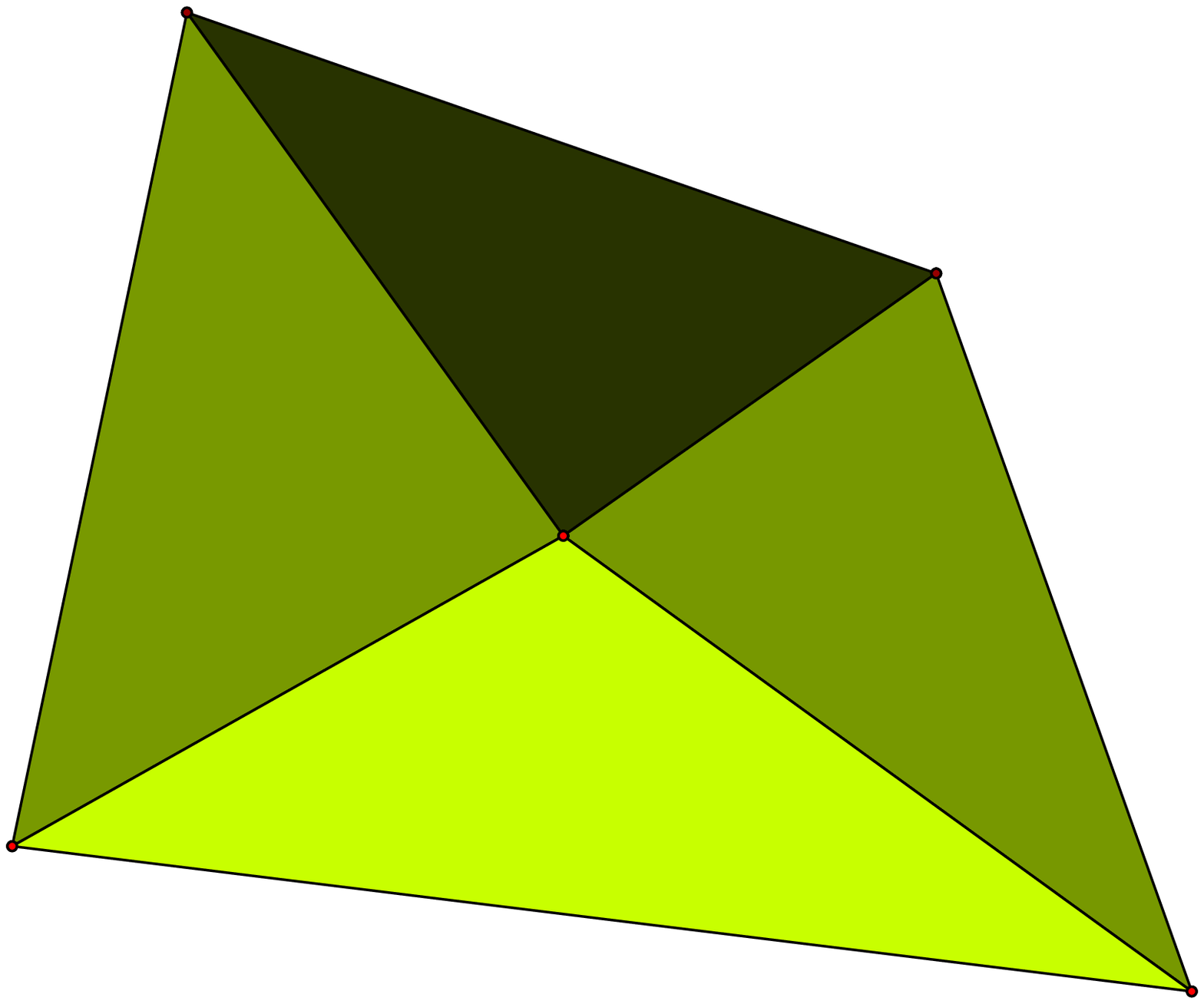,width=1.7in}
		\psfig{figure=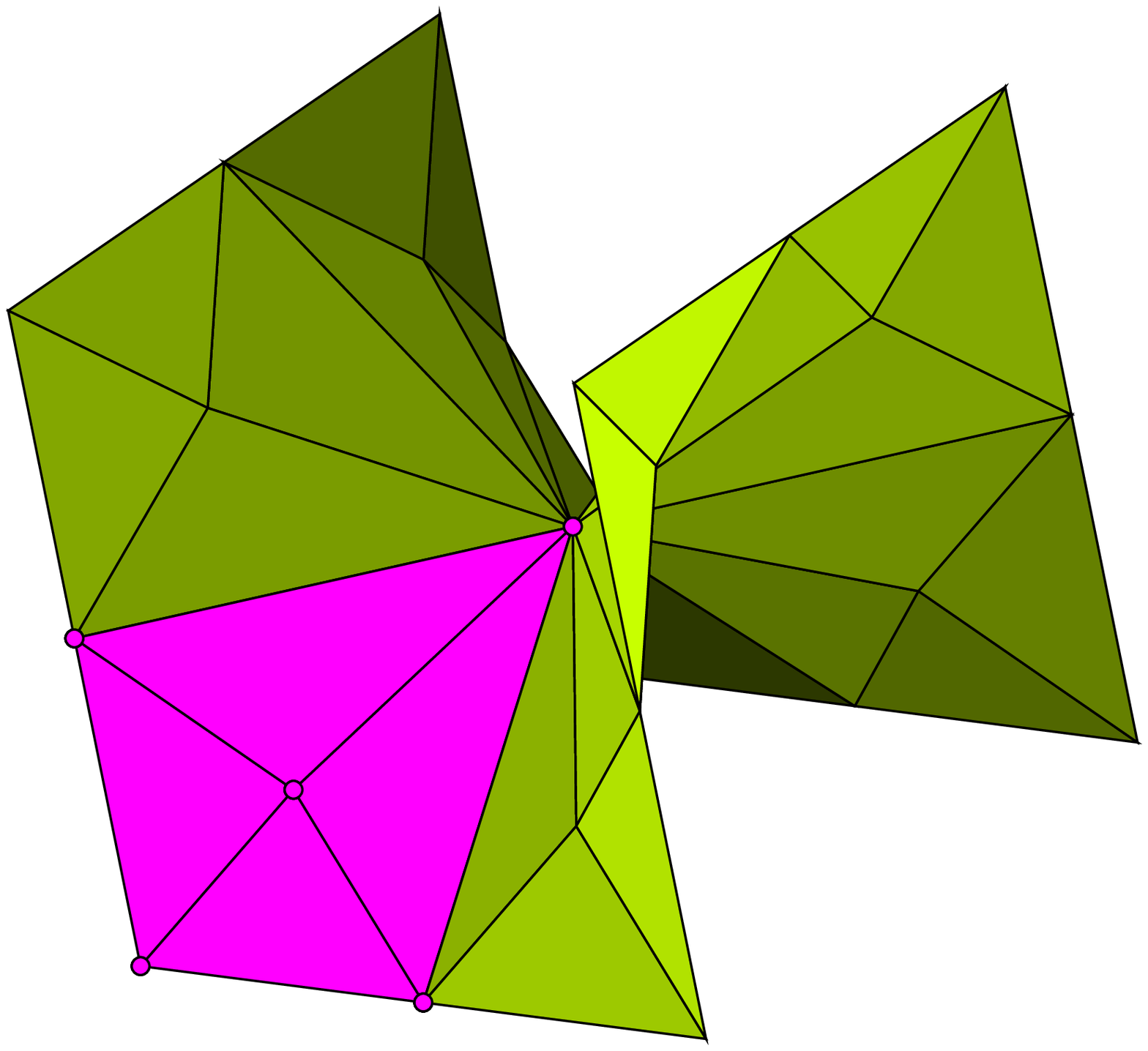,width=1.7in}
	}
}
\centerline{
        \hbox{
		\psfig{figure=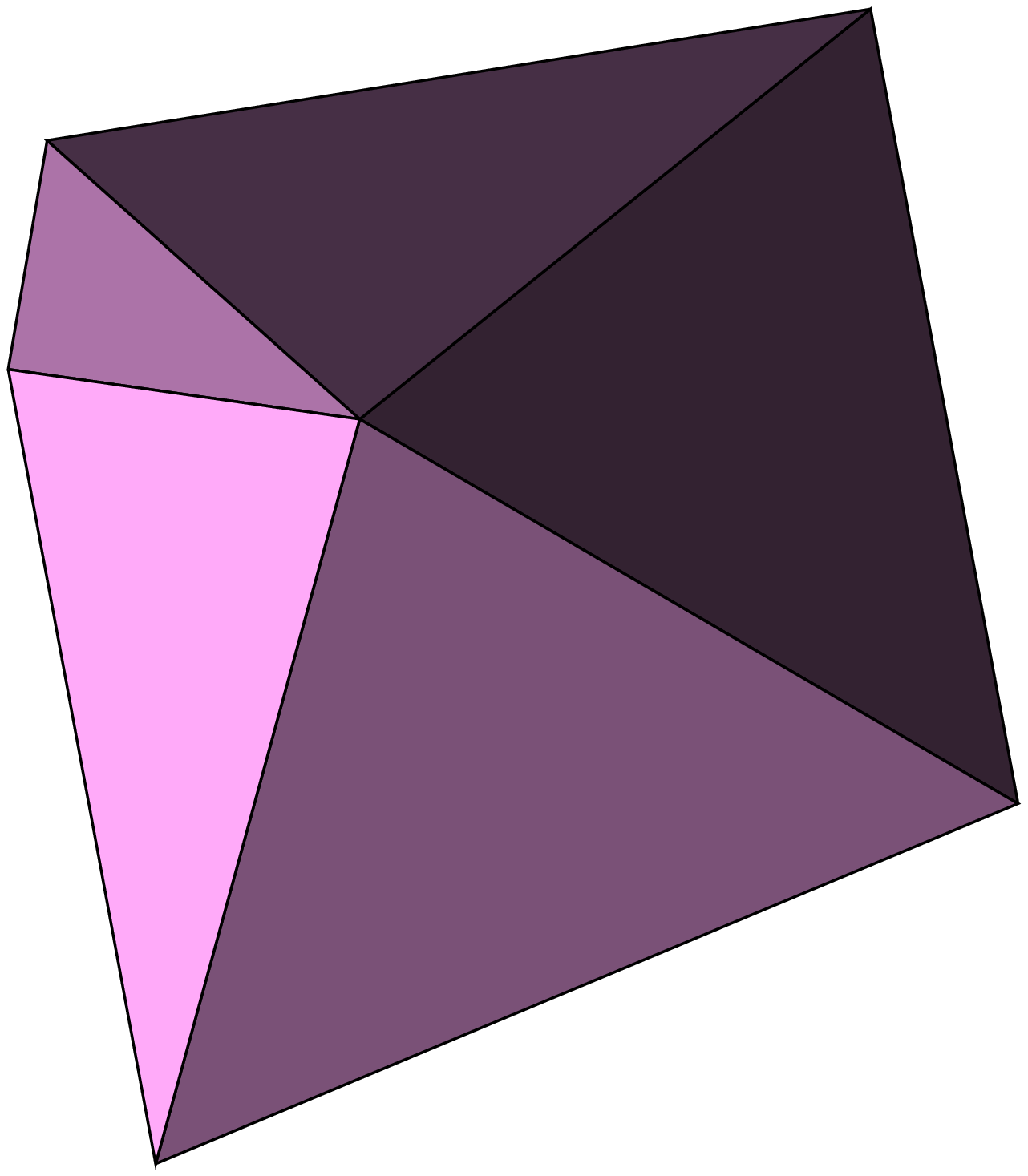,width=1.4in}
		\psfig{figure=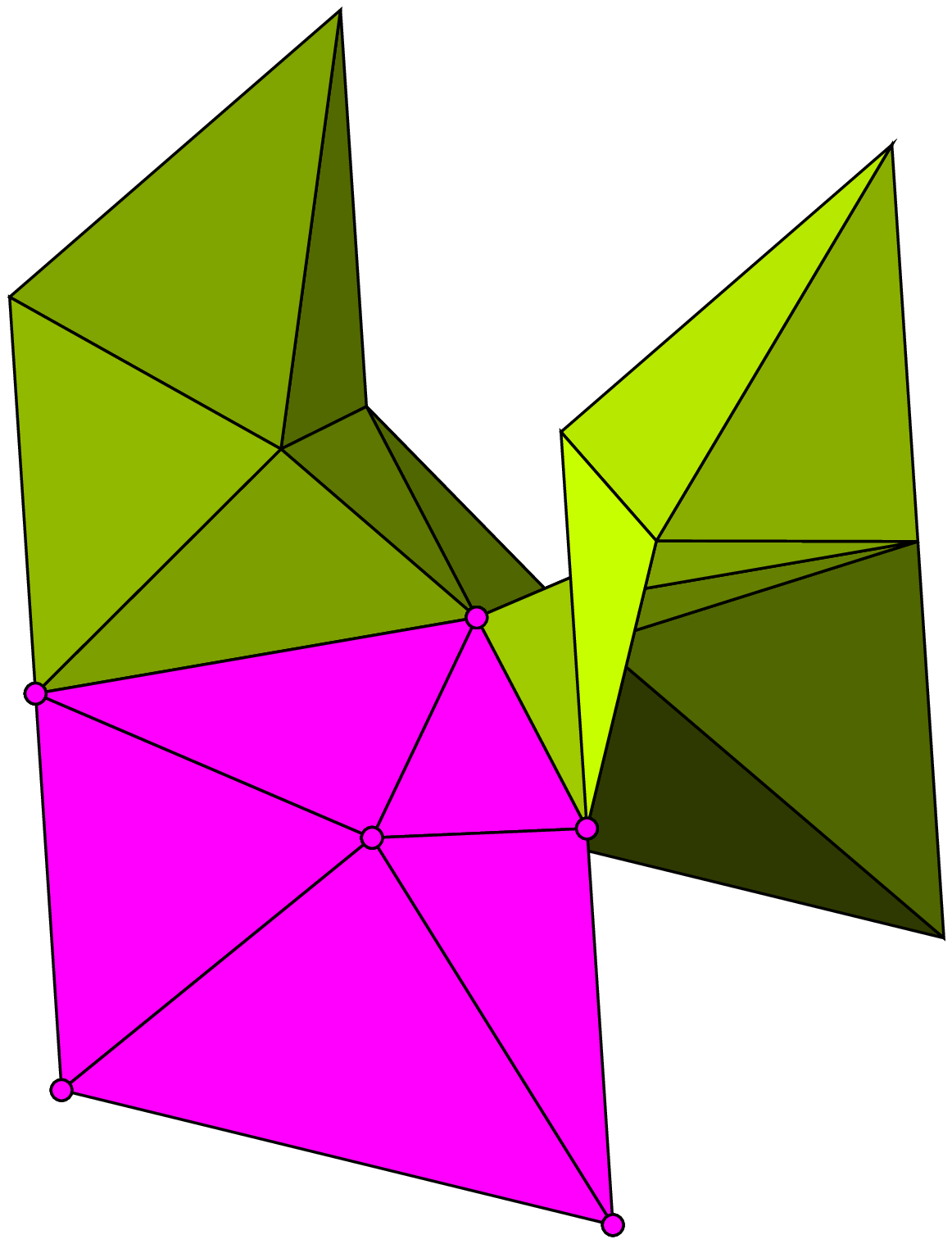,width=1.4in}
		\psfig{figure=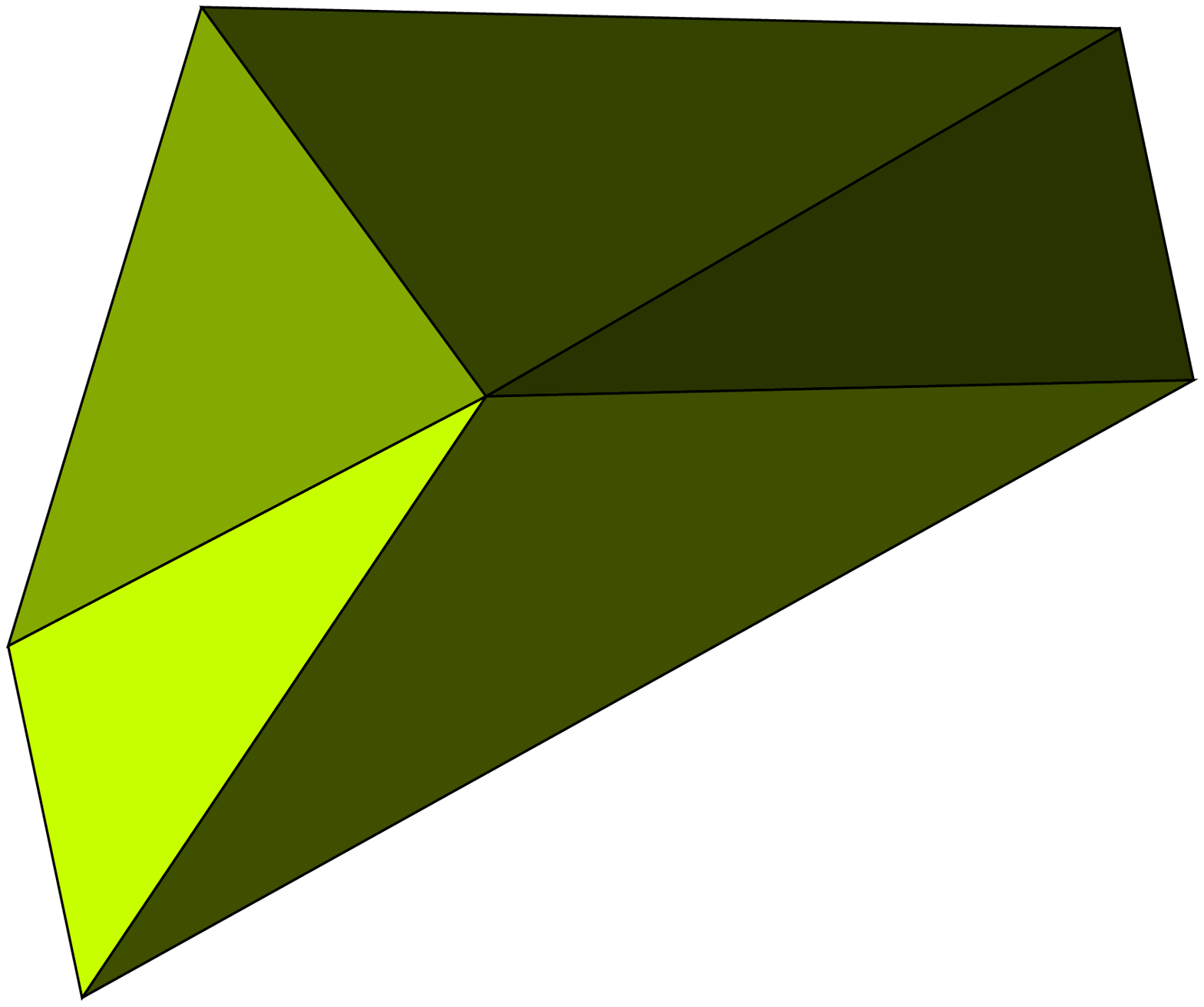,width=1.4in}
	}
}
\caption{\protect\small Four different discrete versions of the superman surface.}
\label{fig:example8}
\end{figure}

The fundamental piece $\mathcal T$ here has eight boundary vertices 
\[ p_1=(1,0,0), \;\; p_2=(1,1,0), \;\; p_3=(1,1,x), \;\; p_4=(0,1,x), \]\[ p_5=(0,1,0), 
\;\; p_6=(0,0,0), \;\; p_7=(0,0,x), \;\; p_8=(1,0,x), \] for any given fixed $x>0$, and 
has one interior vertex \[ p_9=(\tfrac{1}{2},\tfrac{1}{2},\tfrac{x}{2}) \; . \]  
There are eight triangles in $\mathcal T$, which are 
\[ (p_j,p_{j+1},p_9), \; j=1,...,7 \; , \;\;\; (p_8,p_1,p_9). \]
The complete triply-periodic surface is 
generated by including the image of $\mathcal T$ under $180^\circ$ rotations 
about each edge of $\partial {\mathcal T}$, and then continuing to include 
the images under $180^\circ$ rotations about each resulting boundary edge until 
the surface is complete.  It is evident from the symmetries of this surface that 
Equation \eqref{eqn:gradientiszero} holds at every vertex.  
The fundamental piece $\mathcal T$ and a larger part of the resulting complete 
surface are shown on the left-hand side and center of the first row of Figure 
\ref{fig:example8} for $x=1$.  The case for some given $x \in (0,1)$ is shown 
on the right-hand side of the first row of Figure \ref{fig:example8}.  

\subsection{Second example\label{section3-2}} 

The fundamental piece $\mathcal T$ here has six boundary vertices 
\[ p_1=(0,0,0), \;\; p_2=(x,0,0), \;\; p_3=(x,y,0), \]\[ p_4=(x,y,1), 
\;\; p_5=(0,y,1), \;\; p_6=(0,0,1), \] for any given fixed $x,y>0$, and 
has one interior vertex 
\[ p_7=(\tfrac{x}{2},\tfrac{y}{2},\tfrac{1}{2}) \; . \]  There are 
six triangles in $\mathcal T$, which are 
\[ (p_j,p_{j+1},p_7), \; j=1,...,5 \; , \;\;\; (p_6,p_1,p_7). \]
The complete triply-periodic surface is 
generated by $180^\circ$ rotations about boundary edges, just as in the 
previous example.  In this example as well, it is evident from the 
symmetries of this surface that 
Equation \eqref{eqn:gradientiszero} holds at every vertex.  
Two fundamental pieces $\mathcal T$ of different sizes and larger parts of the 
resulting complete surfaces are shown 
in the second row of Figure \ref{fig:example8} ($x=y=1$ in the first case, 
and $x<1<y$ in the second case).  

In the case that $x=y=1$, this fundamental piece $\mathcal T$ has the same 
boundary as a fundamental piece of the smooth Schwarz 
D surface.  Furthermore, for general $x$ and $y$, 
this surface can be viewed as a discrete analog of the superman surface as 
follows: Consider the eight-straight-edged polygonal curve from the point 
$(0,0,-\tfrac{1}{2})$ to the point $(x,-y,-\tfrac{1}{2})$ and then to 
$(x,-y,\tfrac{1}{2})$ and then to $(2 x,0,\tfrac{1}{2})$ 
and then to $(2 x,0,-\tfrac{1}{2})$ and then to $(x,y,-\tfrac{1}{2})$ 
and then to $(x,y,\tfrac{1}{2})$ and then to $(0,0,\tfrac{1}{2})$ 
and then back to $(0,0,-\tfrac{1}{2})$.  This polygonal curve is contained 
in this discrete surface (although not in its edge set) 
and is also the boundary of a smooth superman surface.  

\subsection{Third example\label{section3-3}} 

The fundamental piece $\mathcal T$ here has four boundary vertices 
\[ p_1=(0,0,0), \;\; p_2=(1,1,0), \;\; p_3=(1,1,z), \;\; p_4=(1,0,z), 
\] for any given fixed $z>0$, and has one interior vertex 
\[ p_5=(a,b,c) \; . \]  There are 
four triangles in $\mathcal T$, which are 
\[ (p_j,p_{j+1},p_5), \; j=1,...,3 \; , \;\;\; (p_4,p_1,p_5). \]
Reflecting $\mathcal T$ across the plane $\{ (x_1,0,x_3) \in \mathbb{R}^3 \, | 
\, x_1,x_3 \in \mathbb{R} \}$ and attaching its image to $\mathcal T$, one 
has a larger discrete surface containing eight triangles and six boundary 
edges.  One can extend this larger discrete surface to a complete 
triply-periodic surface by 
$180^\circ$ rotations about boundary edges, just as in the 
previous examples.  In this example, 
Equation \eqref{eqn:gradientiszero} holds at each vertex $p_1$,...,$p_4$ in 
the resulting complete surface.  
However, getting this to hold at $p_5$ requires proper choices of 
$a$ and $b$.  

For simplicity, we restrict to the case $z=1$.  Then, by symmetry, 
we may assume $b=c$.  A computation shows that 
Equation \eqref{eqn:gradientiszero} holding at $p_5$ is equivalent to 
\begin{equation}\label{whatever1} 
(1-a) \sqrt{a^2-2ab+3b^2}=(a-b) \sqrt{(1-a)^2+(1-b)^2} \; , \end{equation}
\begin{equation}\label{whatever2} 
(1-b) \sqrt{a^2-2ab+3b^2}=(3b-a) \sqrt{(1-a)^2+(1-b)^2} \; . \end{equation} 
The solution to this is 
\begin{equation}\label{ab_thirdexample} 
b= \frac{1}{2} \; , \;\;\; a = \frac{3-\sqrt{2}}{2} \; . \end{equation} 
So when $z=1$ and $b=c$ and Equation \eqref{ab_thirdexample} holds, the 
area gradient is zero at each vertex $p_j$ for $j=1,2,...,9$ in the extended 
complete triply-periodic discrete surface, and then symmetries of the surface 
imply the entire complete surface is minimal.  

Since the above minimality condition \eqref{whatever1}-\eqref{whatever2} 
is a system of two equations in two variables $a$ and $b$, we 
say the minimality condition here (when $z=1$) is two-dimensional.  

The fundamental piece $\mathcal T$ with $z=1$ is shown on 
the left-hand side of the third row of Figure \ref{fig:example8}, and a 
larger part of the resulting complete surface is shown just to the right of this.  

\subsection{Fourth example\label{section3-4}} 

The fundamental piece $\mathcal T$ here has five boundary vertices 
\[ p_1=(0,0,0), \;\; p_2=(1,1,0), \;\; p_3=(1,1,z), \;\; p_4=(0,1,z), 
\;\; p_5=(0,0,z) \] for any given fixed $z>0$, and has one interior vertex 
\[ p_6=(a,1-a,b) \; . \]  There are five triangles in $\mathcal T$, which are 
\[ (p_j,p_{j+1},p_6), \; j=1,...,4 \; , \;\;\; (p_5,p_1,p_6). \]
The complete triply-periodic surface is generated by $180^\circ$ rotations 
about boundary edges.  In this example, 
Equation \eqref{eqn:gradientiszero} holds at each vertex $p_1$,...,$p_5$ 
in the resulting complete surface, and 
making it hold also at $p_6$ requires proper choices of $a$ and $b$.  

Like in the previous example, we can 
find a pair of explicit equations, in the variables $a$ and $b$, 
that represent the minimality condition.  These equations are similar to 
those of the previous example, and are slightly more complicated.  One 
can then show the existence of $a$ and $b$ solving this minimality condition.  

Two fundamental pieces $\mathcal T$ of different sizes 
($z=1$ in the first case, and $z<1$ in the second case) are shown on 
the left and right-hand sides of the bottom row of Figure \ref{fig:example8}.  A 
larger part of the resulting complete surface in the case $z=1$ is shown in the 
bottom-middle of Figure \ref{fig:example8}.  

\begin{figure}[tbp]
\centerline{
        \hbox{
		\psfig{figure=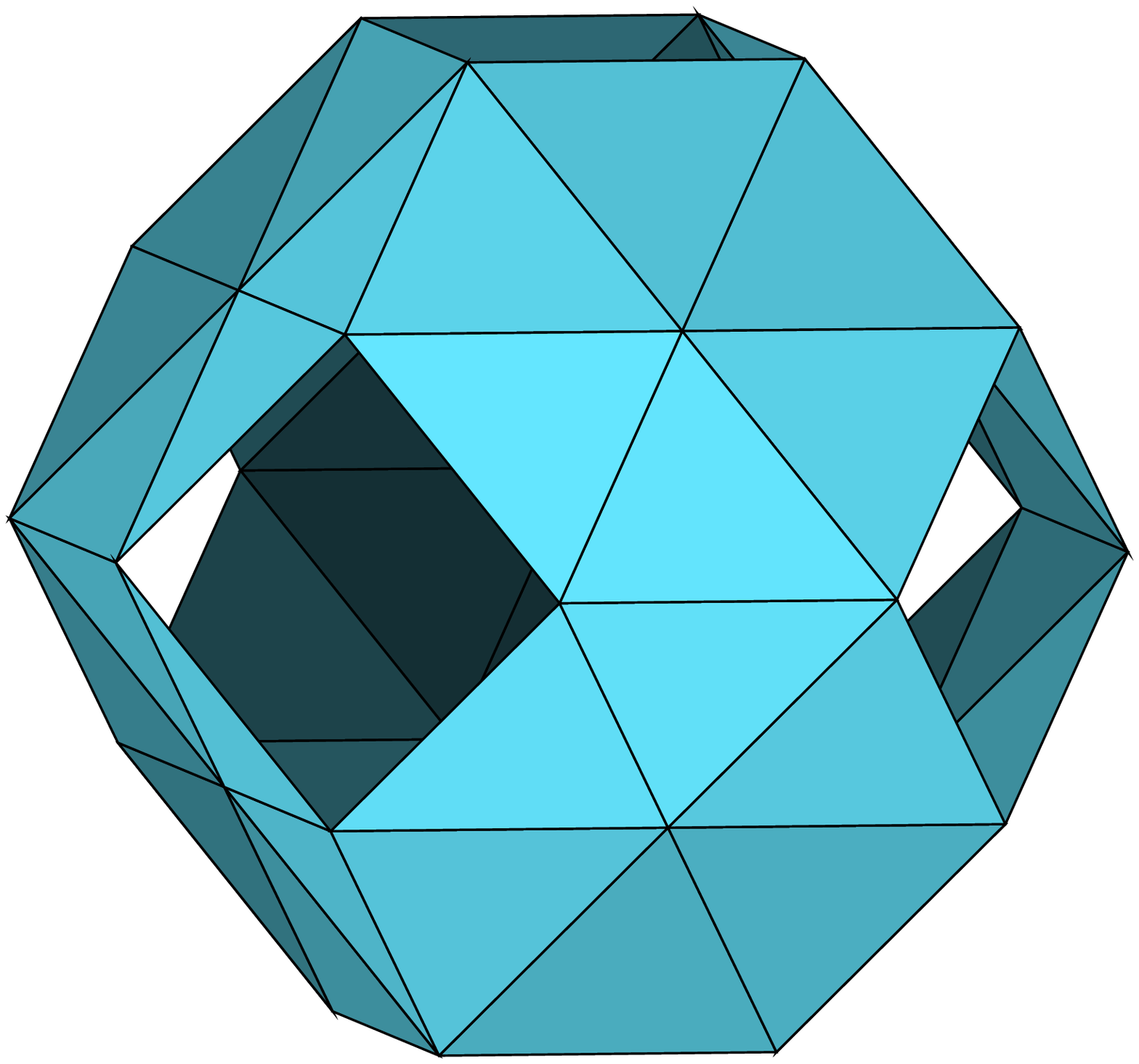,width=1.7in}
		\psfig{figure=exa07x.ps,width=1.7in}
		\psfig{figure=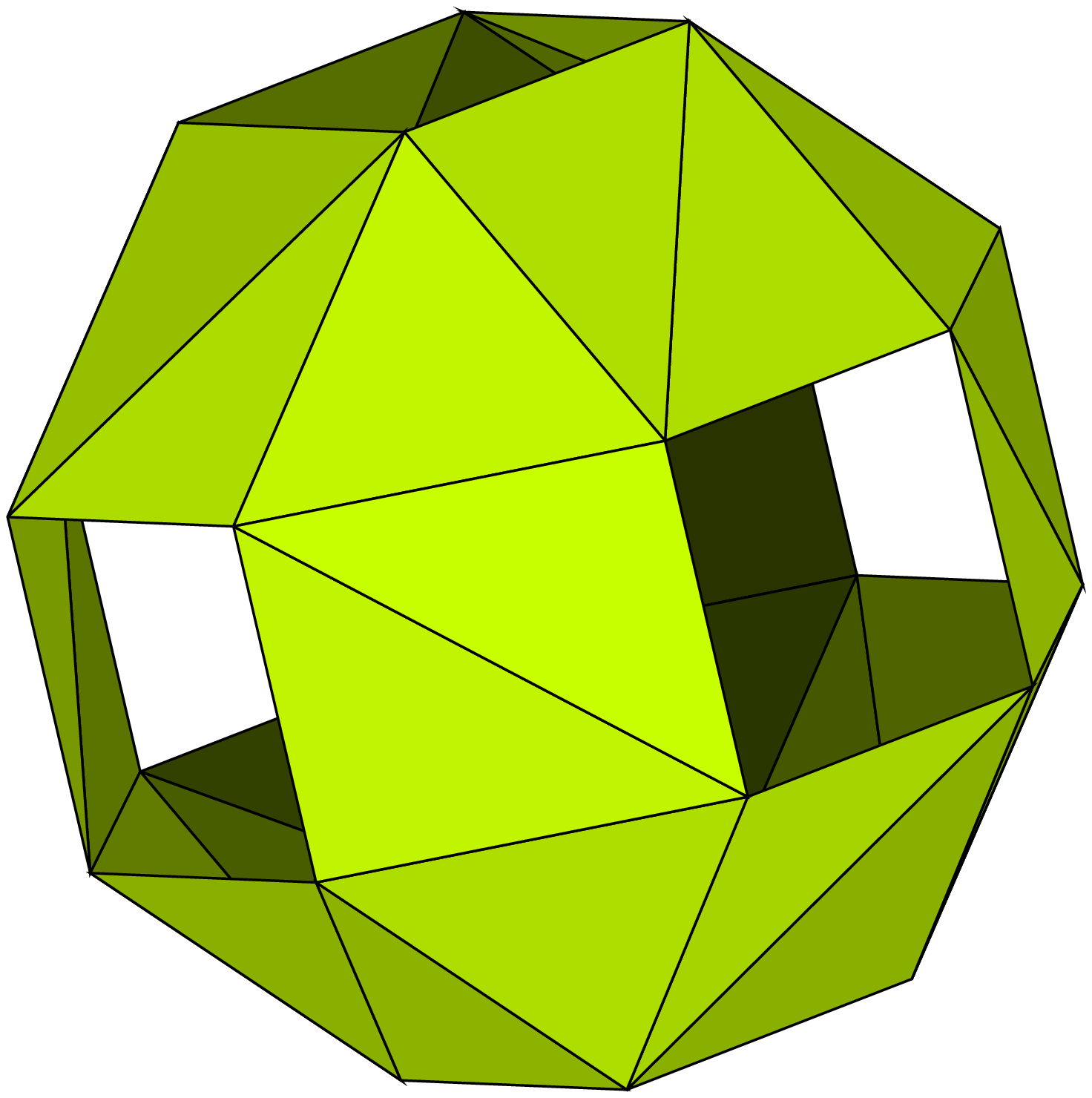,width=1.7in}
	}
}
\centerline{
        \hbox{
		\psfig{figure=exa06y.ps,width=1.3in}
		\psfig{figure=exa07y.ps,width=1.7in}
		\hspace{0.2cm}
		\psfig{figure=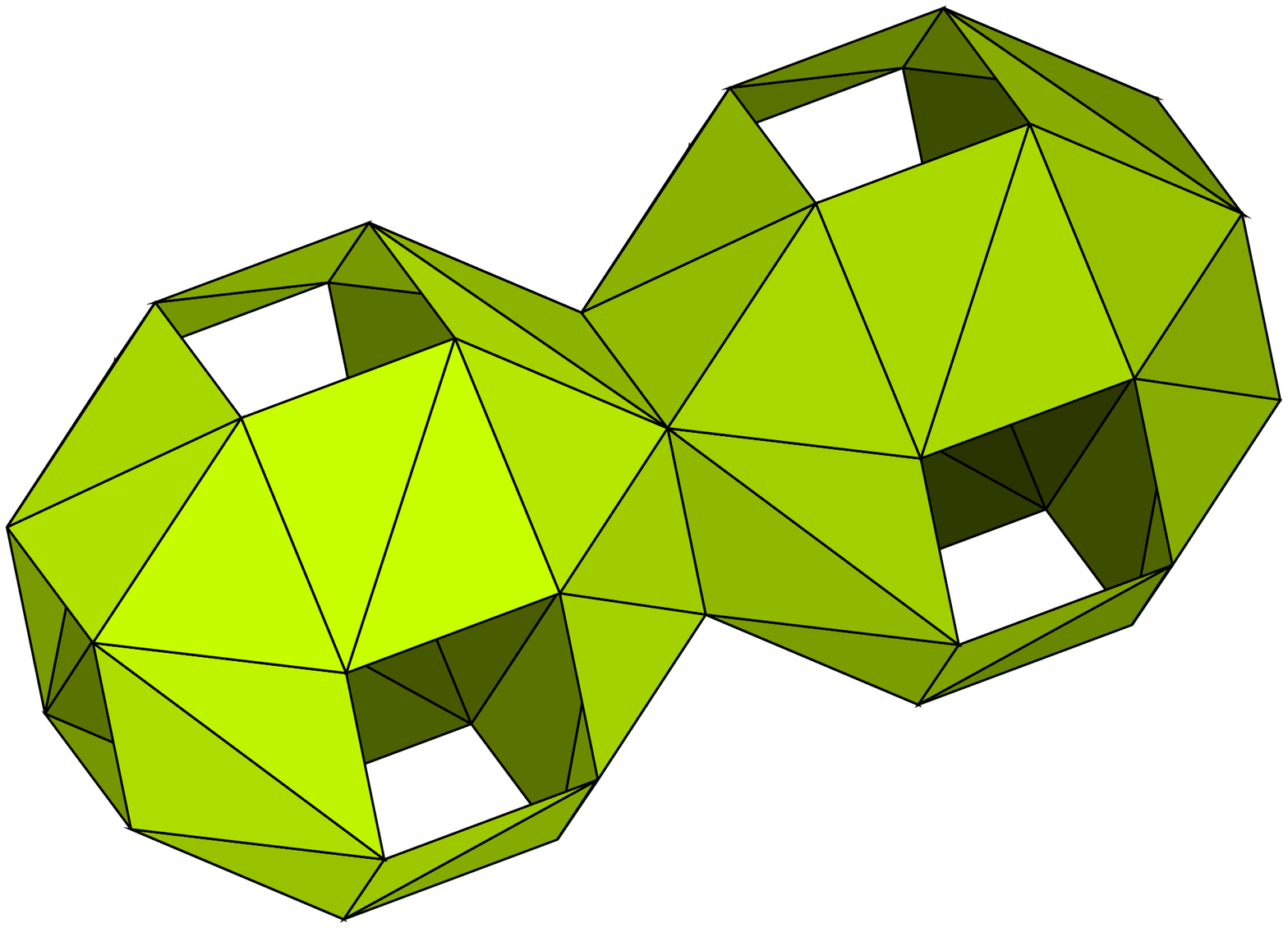,width=1.7in}
	}
}
\caption{\protect\small Three different discrete versions of the Schwarz P surface.}
\label{fig:example6}
\end{figure}

\section{Discrete versions of the Schwarz P surface\label{section4}}

\subsection{First two examples\label{section4-1}} 

Consider the vertices 
\[ p_1=(3,0,6), \;\; p_2=(6,0,3), \;\; p_3=(6,3,0), \]\[ p_4=(3,6,0), 
\;\; p_5=(0,6,3), \;\; p_6=(0,3,6), \;\; p_7=(3,3,3), \] and let 
${\mathcal T}_1$ be the planar fundamental domain with the six triangles 
\[ (p_j,p_{j+1},p_7), \; j=1,...,5 \; , \;\;\; (p_6,p_1,p_7). \]
Also, consider the vertices 
\[ p_1=(3,0,6), \;\; p_2=(4,0,4), \;\; p_3=(6,0,3), \;\; p_4=(6,2,2), \]\[ 
p_5=(6,3,0), \;\; p_6=(4,4,0), \;\; p_7=(3,6,0), \;\; p_8=(2,6,2), \]\[ 
p_9=(0,6,3), \;\; p_{10}=(0,4,4), \;\; p_{11}=(0,3,6), 
\;\; p_{12}=(2,2,6), \;\; p_{13}=(3,3,3), \] and let 
${\mathcal T}_2$ be the fundamental domain with the twelve triangles 
\[ (p_j,p_{j+1},p_{13}), \; j=1,...,11 \; , \;\;\; (p_{12},p_1,p_{13}). \]

We can extend ${\mathcal T}_j$ (for either $j=1,2$) to a complete discrete 
surface by including the images of ${\mathcal T}_j$ under the reflections 
across the planes $\{(x,y,6k) \, | \, x,y \in {\mathbb{R}} \}$, 
$\{(x,6k,z) \, | \, x,z \in {\mathbb{R}} \}$ and 
$\{(6k,y,z) \, | \, y,z \in {\mathbb{R}} \}$ for 
all integers $k$.  Furthermore, Equation \eqref{eqn:gradientiszero} holds 
at all vertices of the extended surface, so it is minimal.  (See the first two 
columns of Figure \ref{fig:example6}.)  

The surface produced by ${\mathcal T}_1$ (resp. ${\mathcal T}_2$) is a simpler 
(resp. more complicated) version of a discrete 
Schwarz surface.  Note that they are analogous to the bottom-right picture in 
Figure \ref{fig:example11}.  (The second 
example ${\mathcal T}_2$ was also shown in \cite{PR}.)

\subsection{Third example\label{section4-2}}

Consider the ten vertices 
\[ p_j=(a,(-1)^{j} a,1), \;\; p_{j+2}=(a,1,(-1)^{j+1} a), 
\;\; p_{j+4}=(a,(-1)^{j+1} a,-1), \]\[ p_{j+6}=(1,(-1)^{j} a,a), \;\; 
p_{j+8}=(1,(-1)^{j+1} a,-a), 
\;\;\;\;\;\;\;\;\;\;\;\;\;\;\;\;\;\;\;\;\;\;\;\;\; j=1,2, \] and let 
$\hat{\mathcal T}$ be the discrete surface with the eight triangles 
\[ (p_1,p_2,p_7), \;\; (p_2,p_8,p_7), \;\; (p_2,p_3,p_8), \;\; (p_3,p_4,p_8), \]
\[ (p_4,p_9,p_8), \;\; (p_4,p_5,p_9), \;\; (p_5,p_6,p_9), \;\; (p_6,p_{10},p_9). \]
Then let $\mathcal T$ be the discrete surface, with $24$ vertices and $32$ 
faces, that is made by including the four images of $\hat{\mathcal T}$ under 
the rotations about the axis $\{(0,0,r) \, | \, r \in \mathbb{R} \}$ of angles 
$0^\circ$, $90^\circ$, $180^\circ$ and $270^\circ$.  
This $\mathcal T$ is shown in the upper-right of Figure \ref{fig:example6}.  

One can then generate a complete triply-periodic surface by including the 
images of $\mathcal T$ under reflections across the planes 
$\{(x,y,k) \, | \, x,y \in {\mathbb{R}} \}$, 
$\{(x,k,z) \, | \, x,z \in {\mathbb{R}} \}$ and 
$\{(k,y,z) \, | \, y,z \in {\mathbb{R}} \}$ for 
all integers $k$.  The result of applying one such reflection is shown in 
the lower picture of the right-most column of Figure \ref{fig:example6}.  

The condition for this discrete triply-periodic surface to be minimal is that 
\[ a= \frac{3 \sqrt{2}-\sqrt{3}}{6 \sqrt{2}-\sqrt{3}} \; , 
\]  i.e. for this value of $a$, Equation 
\eqref{eqn:gradientiszero} holds at every vertex of the surface.  

\subsection{Examples based on discrete minimal catenoids\label{section4-3}}

Here we give two closely-related types of examples based on discrete minimal 
catenoids.  One type is a discrete analog of the Schwarz P surface.  The 
other type is actually an analog of the smooth Schwarz H surface, not the 
Schwarz P surface.  
To construct these examples, we will use discrete analogs of the catenoid 
\cite{PR}, which are 
described in terms of the hyperbolic cosine function, just as the smooth 
catenoid was in Equation \eqref{smoothcatenoid}.  

The vertices of a discrete minimal catenoid lie on congruent
planar polygonal meridians, and the meridians are contained in planes that 
meet along a single line (the axis) at equal angles.  Every meridian is the image 
of every other meridian by some rotation about the axis.  
By drawing edges between corresponding vertices of adjacent meridians (i.e. so 
that these edges are perpendicular to the axis), we have 
a piecewise linear continuous surface tessellated by planar isosceles 
trapezoids.  We can 
triangulate each trapezoid any way we please without affecting minimality, 
as noted in Remark \ref{flatpartsofminimalsurfaces}, so we shall triangulate 
each trapezoid by drawing a single diagonal edge across it.  

\begin{figure}[tbp]
\centerline{
        \hbox{
		\psfig{figure=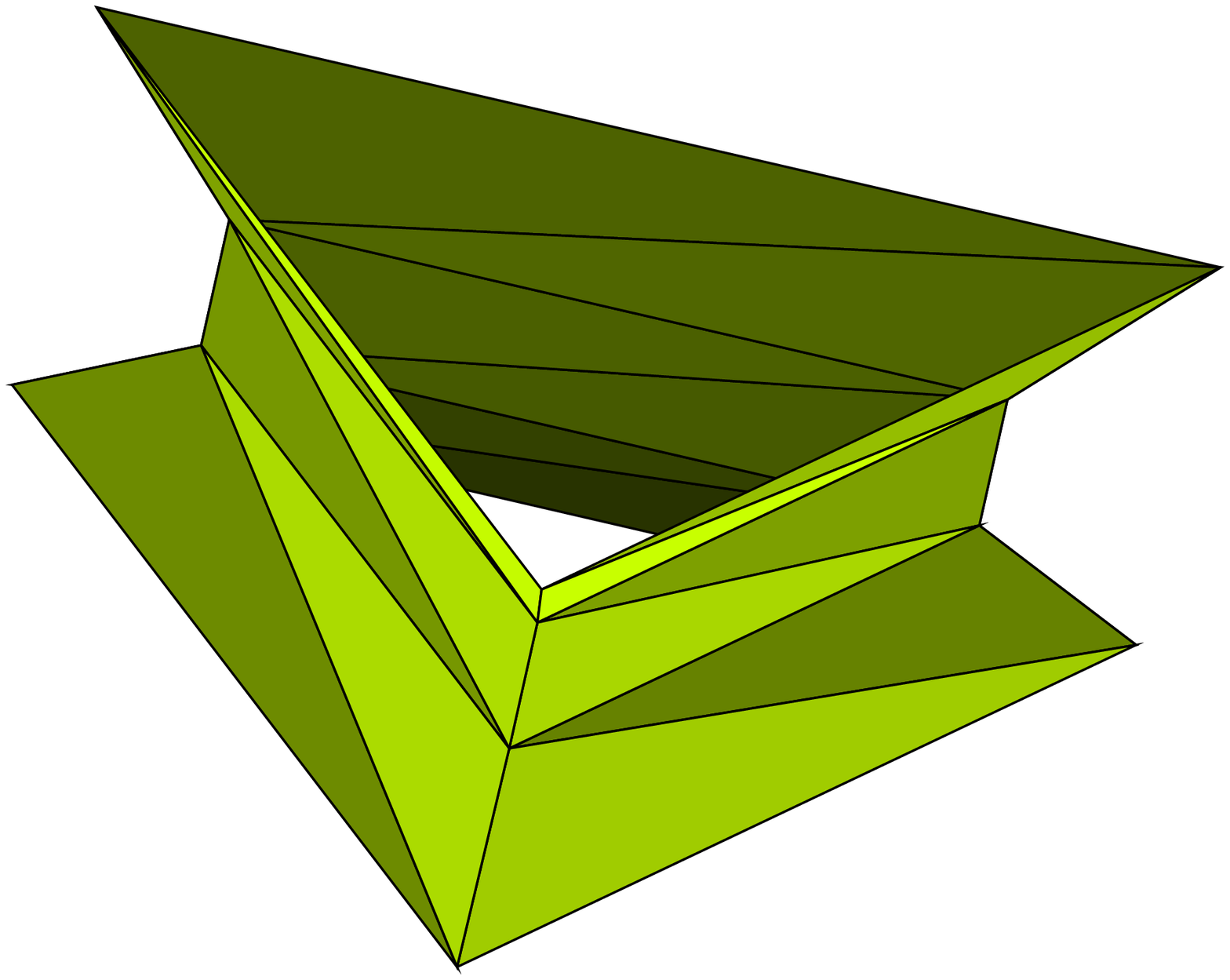,width=1.4in}
		\psfig{figure=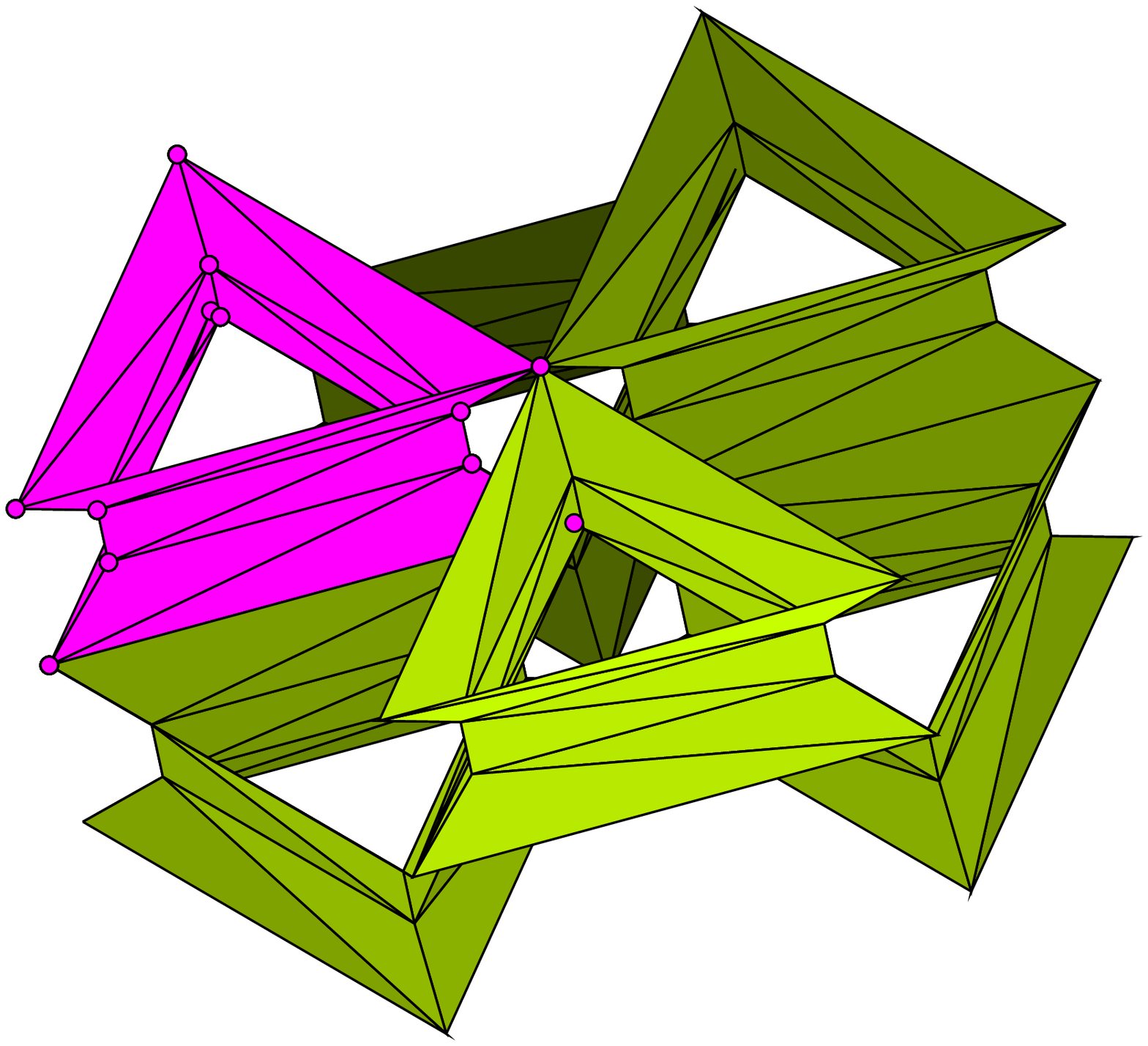,width=1.6in}
		\psfig{figure=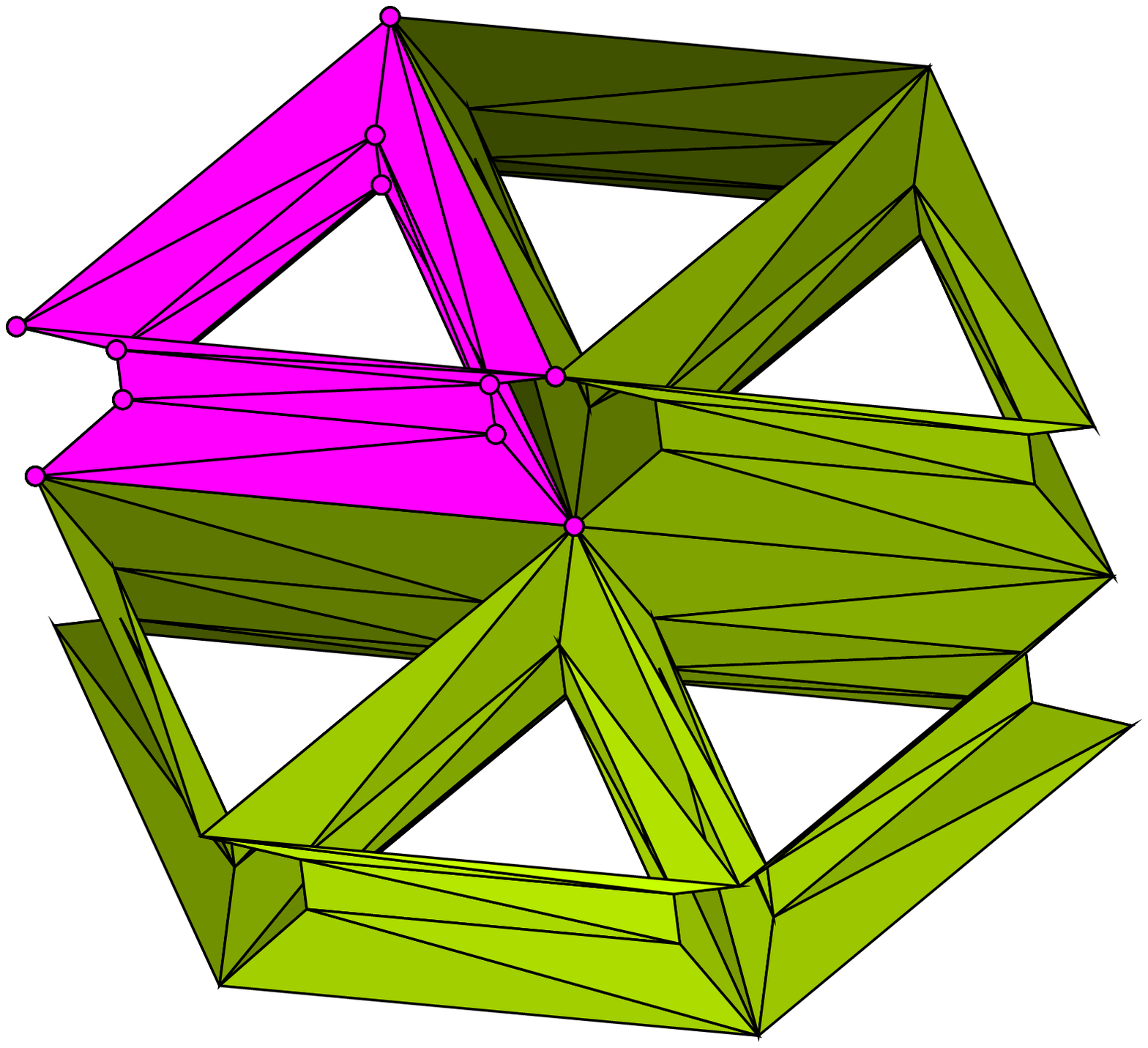,width=1.6in}
	}
}
\caption{\protect\small Discrete version of the Schwarz H surface.}
\label{fig:example11b}
\end{figure}

Two examples of discrete catenoids are shown in the first two pictures in the 
upper row of Figure \ref{fig:example11}.  Both of these pictures have adjacent 
meridians in planes meeting at $90^\circ$ angles.  The first (resp. second) one 
has four (resp. five) vertices in each meridian.  Another example is shown in 
the left-most picture of Figure \ref{fig:example11b}, where the adjacent 
meridians lie in planes meeting at $120^\circ$ angles, and there are four 
vertices in each meridian.  

To explicitly describe discrete catenoids, we need only specify: 
\begin{enumerate}
\item The axis $\ell$: let us fix $\ell = \{(0,0,z) \, | \, z \in 
      {\mathbb{R}} \}$.  
\item The angle $\theta$ between planes of adjacent meridians: let us fix
      $\theta = \frac{2 \pi}{k}$ for some integer $k \geq 3$.  
\item The locations of the vertices along one meridian.  
\end{enumerate}
We can place one meridian in the plane 
$\{(x,0,z) \, | \, x,z \in {\mathbb{R}} \}$, and locating its vertices 
at the following points will ensure minimality of the surface (see \cite{PR}):  
\[ p_j=(r\cosh \left( {\frac{1}{r}}a (z_{0}+j\delta) \right),0,z_{0}+j\delta)
\] with $j=j_0,j_0+1,...,j_1$ for some integers $j_0$ and $j_1$ ($j_0<j_1$), and 
with 
\[
a={\frac{r}{\delta }}\func{arccosh}\left( 1+{\frac{1}{r^{2}}}{\frac{\delta
^{2}}{1+\cos \theta }}\right) , 
\]
where $r>0$ and $\delta >0$ and $z_0 \in \mathbb{R}$ are constant.  
The edges along this meridian are $\overline{p_j p_{j+1}}$ 
for $j$ between $j_0$ and $j_1-1$.  

For our application, we shall restrict to either $k=4$, as in Figure 
\ref{fig:example11}, or to $k=3$, as in Figure \ref{fig:example11b}.  
We shall further assume that either 
\begin{itemize}
\item $z_0=0$ and $j_0 = -j_1 < 0$, or 
\item $z_0=\tfrac{\delta}{2}$ and $j_0 = -j_1-1 < -1$.  
\end{itemize}
Either of these conditions will produce a discrete minimal surface $\mathcal T$ 
whose trace has dihedral symmetry.  One can then extend $\mathcal T$ by 
$180^\circ$ rotation about boundary lines to a complete embedded discrete surface 
in $\mathbb{R}^3$.  To conclude minimality of this complete surface, it remains 
only to check that Equation \eqref{eqn:gradientiszero} holds at any vertex 
contained in any edge about which a $180^\circ$ rotation was made, and this is 
clear from the symmetry of the surface.  

The case when $k=4$ and $z_0=0$ and $j_0=-j_1=-2$ is shown in the second 
picture of the first row of Figure \ref{fig:example11}, and a larger portion 
of the resulting complete minimal surface is shown in the picture just to the 
right of it.  The case when $k=3$ and $z_0=\tfrac{\delta}{2}$ and $j_0=-j_1-1=-2$ 
is shown in the left-most 
picture of Figure \ref{fig:example11b}, and a larger portion 
of the resulting complete minimal surface is shown in the middle of Figure 
\ref{fig:example11b}.  When $k=4$, the analogy to the smooth Schwarz P surface 
is clear.  When $k=3$, one can imagine a smooth embedded minimal annulus 
with the same boundary as $\mathcal T$, and this surface is called 
the Schwarz H surface.  

As explained in Section \ref{section1}, there are infinitely many different ways 
(by using combinations of reflections and $180^\circ$ rotations that are 
not allowed in the smooth case) to extend 
$\mathcal T$ to a complete discrete minimal surface.  Two such ways are shown 
in the upper right of Figure \ref{fig:example11}, and another two 
ways are shown in the center and right-hand side of 
Figure \ref{fig:example11b}.  The two examples in Figure \ref{fig:example11} and 
the central one in Figure \ref{fig:example11b} can be extended to complete 
triply-periodic discrete minimal surfaces by $180^\circ$ rotations about 
boundary edges.  The right-most example in Figure \ref{fig:example11b} can be 
extended to a complete triply-periodic discrete minimal surface by using 
horizontal translations perpendicular to the axis $\ell$ that generate a 
$2$-dimensional hexagonal grid, and then by applying vertical 
translations parallel to $\ell$ of length $2 \delta (j_1-j_0)$.  The upper-right 
examples in both Figures \ref{fig:example11} and \ref{fig:example11b} 
are applications of {\bf Method 2}.  

\begin{remark}
When $k=4$ and $j_0=-j_1=1$, and when $r$ and $\delta$ are chosen properly, this 
$\mathcal T$ can produce the same surface as ${\mathcal T}_1$ produced in 
Subsection \ref{section4-1}.  The way of triangulating the planar isosceles 
trapezoids was different in Subsection \ref{section4-1}, but by Remark 
\ref{flatpartsofminimalsurfaces} this is irrelevant to the minimality of 
the surfaces, and the two examples are the same in the sense that they have the 
same traces in ${\mathbb{R}^{3}}$.  
\end{remark}

\begin{figure}[tbp]
\centerline{
        \hbox{
		\psfig{figure=exa01a.ps,width=1.7in}
		\psfig{figure=exa01b.ps,width=1.7in}
		\psfig{figure=exa01c.ps,width=1.7in}
	}
}
\centerline{
        \hbox{
		\psfig{figure=exa01aa.ps,width=1.7in}
		\psfig{figure=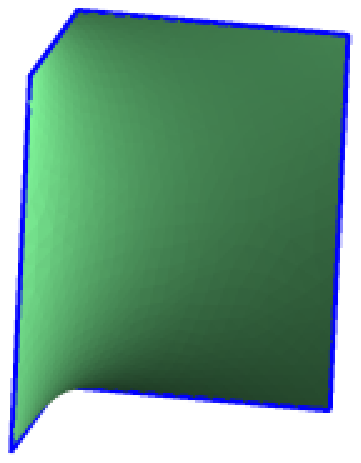,width=1.7in}
		\psfig{figure=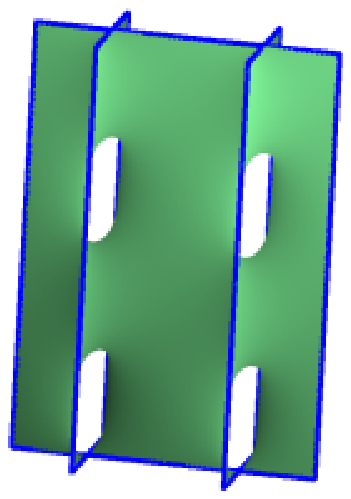,width=1.7in}
	}
}
\caption{\protect\small Discrete and smooth Schwarz CLP surfaces.}
\label{fig:example1}
\end{figure}

\section{Other examples\label{section5}}

\subsection{Discrete Schwarz CLP surface\label{section5-1}} 

The fundamental piece $\mathcal T$ here has six boundary vertices 
\[ p_1=(x,0,0), \;\; p_2=(0,0,0), \;\; p_3=(0,y,0), \]\[ p_4=(0,y,1), \;\; p_5=(0,0,1), 
\;\; p_6=(x,0,1) \] for any given fixed $x,y>0$, and has one interior vertex 
\[ p_7=(a,b,\tfrac{1}{2}) \; . \]  There are six triangles in $\mathcal T$, which are 
\[ (p_j,p_{j+1},p_7), \; j=1,...,5 \; , \;\;\; (p_6,p_1,p_7) \; . \]
The complete triply-periodic surface is generated by $180^\circ$ rotations 
about boundary edges, continuing to make such rotations until the 
surface is complete.  Every vertex in $\partial {\mathcal T}$, and every vertex 
that is an image of a vertex in $\partial {\mathcal T}$ under these rotations, 
satisfies Equation \eqref{eqn:gradientiszero}, because of the symmetry of the 
surface.  The condition for Equation \eqref{eqn:gradientiszero} to hold at 
the interior vertex $p_7$ and all images of $p_7$ under these rotations is that 
\[\frac{2 y a}{\sqrt{a^2+\frac{1}{4}}} + 
\frac{a}{\sqrt{a^2+(y-b)^2}} + \frac{a-x}{\sqrt{b^2+(x-a)^2}} = 0 \; , \] 
\[\frac{2 x b}{\sqrt{b^2+\frac{1}{4}}} + 
\frac{b}{\sqrt{b^2+(x-a)^2}} + \frac{b-y}{\sqrt{a^2+(y-b)^2}} = 0 \; . \]
When $x=y=\tfrac{\sqrt{2}}{2}$, one explicit solution is $a=b=\frac{\sqrt{2}-1}{2}$.  
The fundamental piece $\mathcal T$ and a larger part of the resulting complete 
surface are shown in the left-most column of 
Figure \ref{fig:example1} for these values of $x$, $y$, $a$, and $b$.  

\begin{figure}[tbp]
\centerline{
        \hbox{
		\psfig{figure=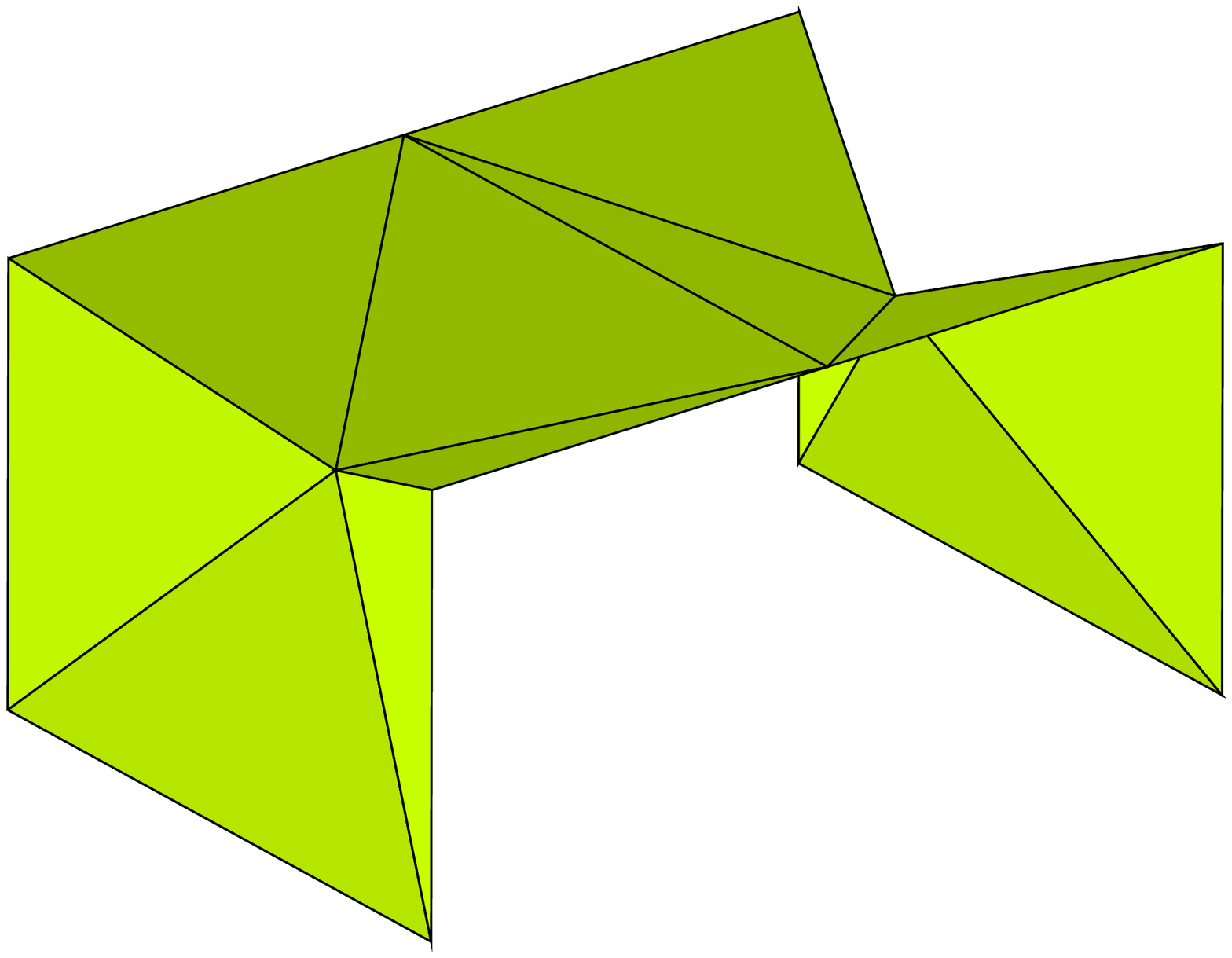,width=1.7in}
		\psfig{figure=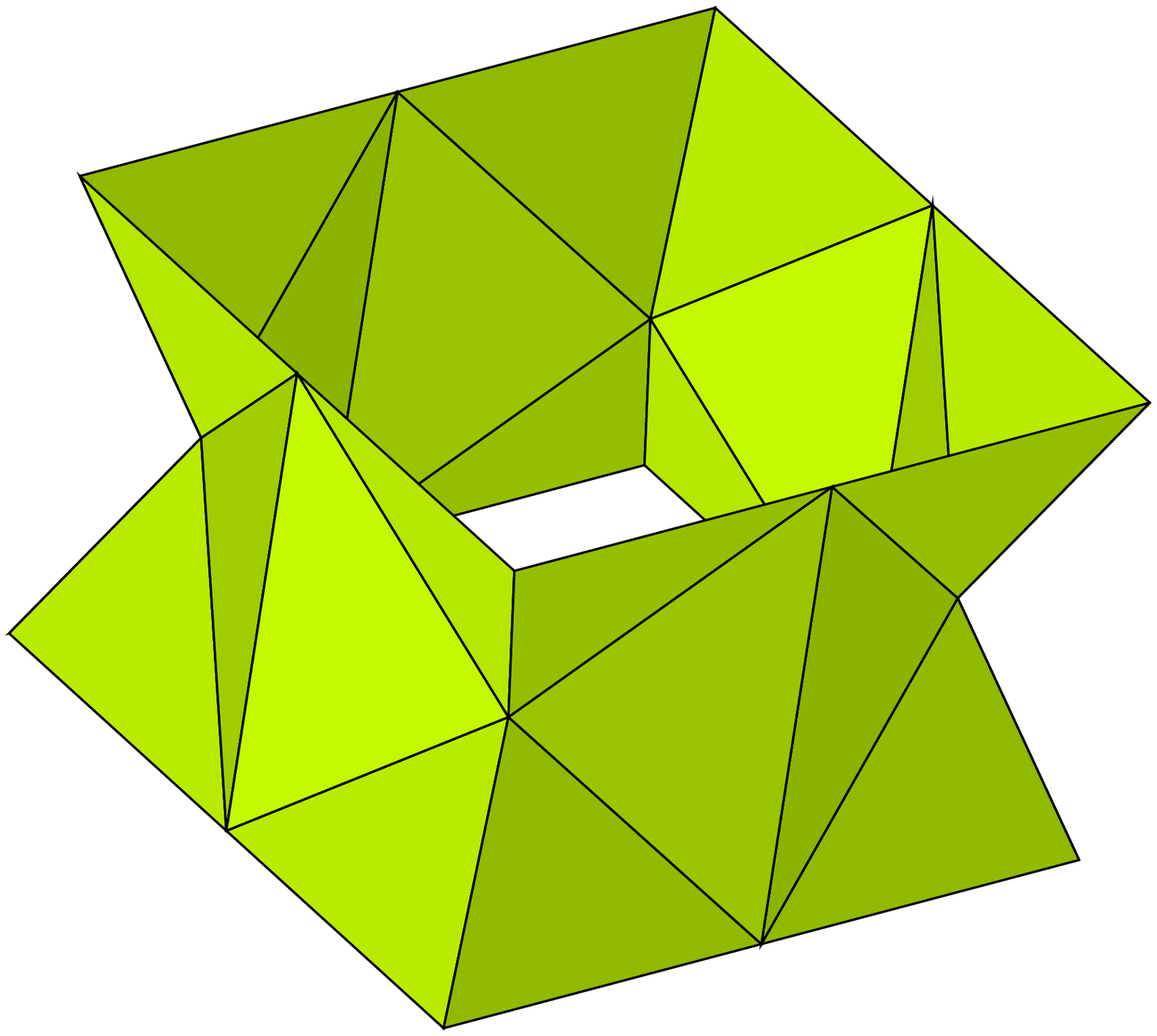,width=1.7in}
	}
}
\caption{\protect\small Variants of the discrete Schwarz CLP surface.}
\label{fig:example1a}
\end{figure}

For general choices of $x$ and $y$, there is always a solution to the above 
system of two equations with respect to the two variables $a$ and $b$, 
thus giving $\nabla_{p_7} \func{area} {\mathcal T}=0$.  Thus, for 
general $x$ and $y$, the minimality condition for this example is 
two-dimensional.  Fundamental pieces $\mathcal T$ for other 
choices of $x$ and $y$ are shown in the upper-center and upper-right of Figure 
\ref{fig:example1} ($x=y$ in the center and $x \neq y$ on the right).  

We can also apply {\bf Method 2} here.  
For example, suppose we include the reflection of $\mathcal T$ across the plane 
$P$ containing the three points $(x,0,0)$, $(x,1,0)$, $(x,0,1)$ 
along with $\mathcal T$ to get a discrete minimal surface ${\mathcal T}_1$ with 
twelve triangles, see the left-hand side of Figure 
\ref{fig:example1a}.  (Such a reflection across $P$ would not be allowed for 
the smooth Schwarz CLP surface.)  We can then extend $\mathcal{T}_1$ 
to a complete triply-periodic discrete minimal surface by $180^\circ$ 
rotations about boundary edges, and this surface is yet another 
discrete superman surface.  

For a second example of applying {\bf Method 2}, 
suppose we include the reflection of ${\mathcal T}_1$ across the plane 
$Q$ containing the three points $(0,y,0)$, $(1,y,0)$, $(0,y,1)$ 
along with ${\mathcal T}_1$ to get a discrete minimal surface ${\mathcal T}_2$ with 
twenty-four triangles, see the right-hand side of Figure 
\ref{fig:example1a}.  (Such a reflection again would not be allowed in 
the smooth case.)  We can then extend $\mathcal{T}_2$ 
to a complete triply-periodic discrete minimal surface by $180^\circ$ 
rotations about boundary edges, and this gives yet another 
discrete Schwarz P surface.  

\subsection{Discrete I-Wp and F-Rd surfaces\label{section5-2}} 

\begin{figure}[tbp]
\centerline{
        \hbox{
		\psfig{figure=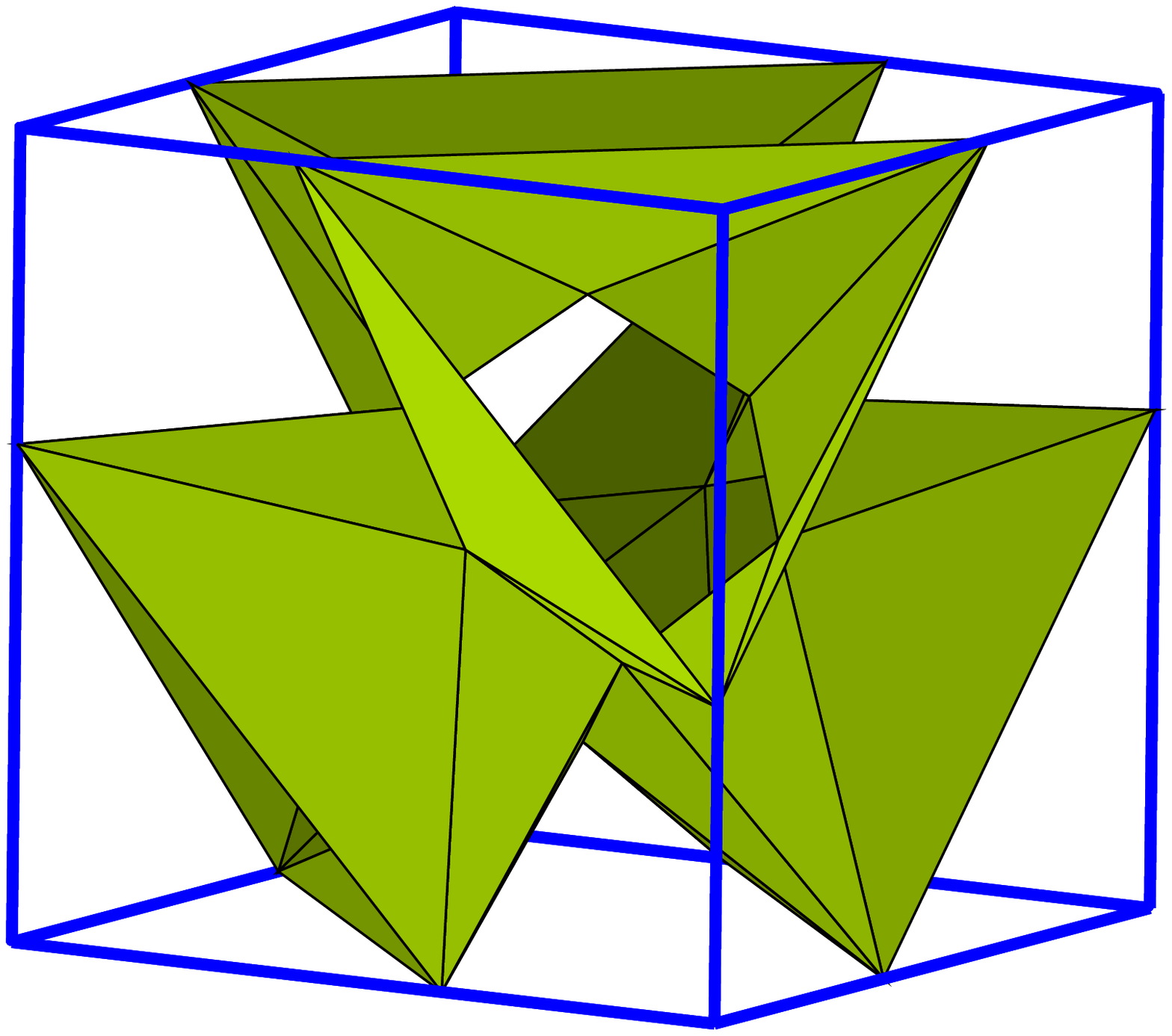,width=2.1in}
		\psfig{figure=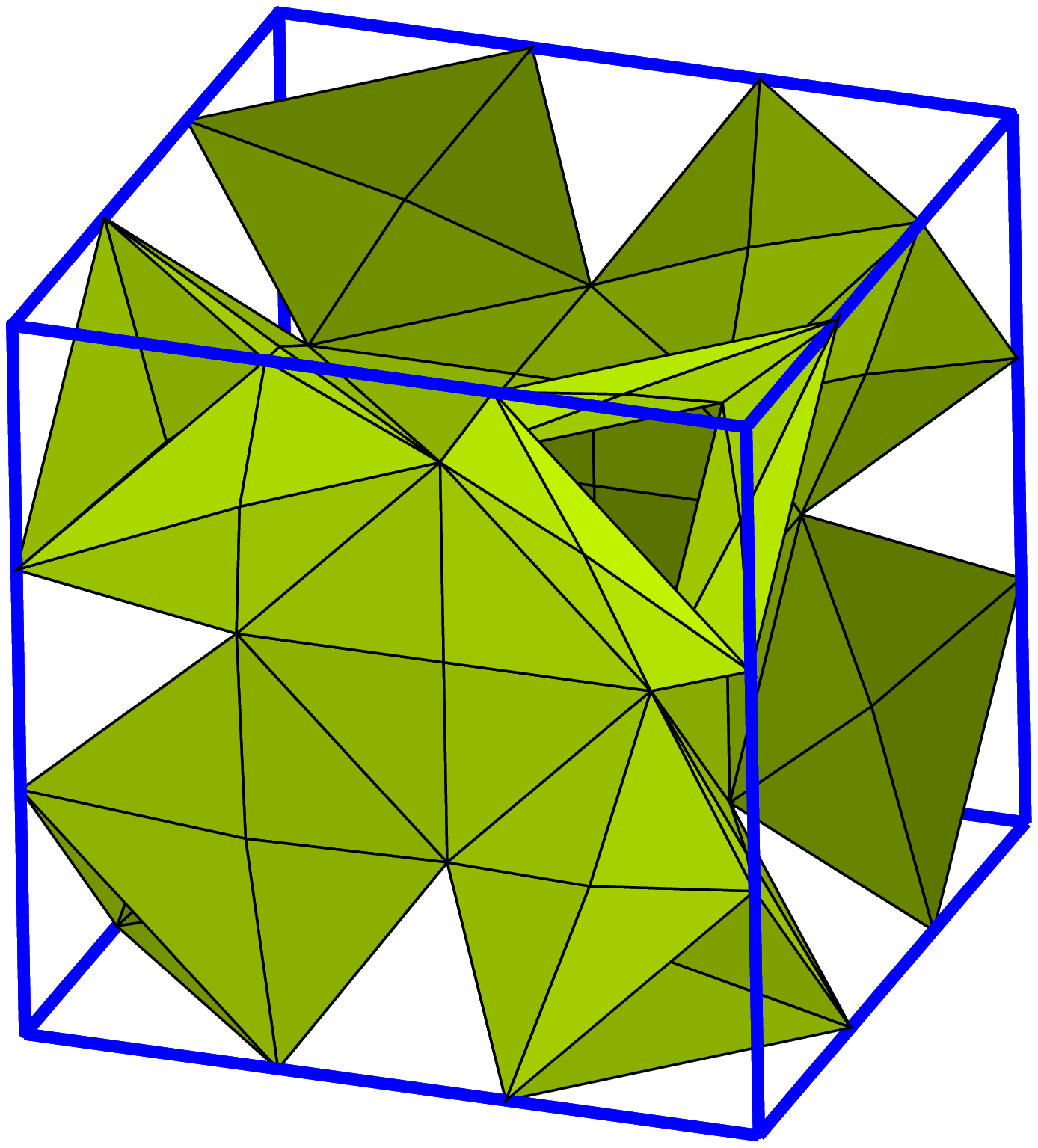,width=1.7in}
	}
}
\caption{\protect\small Discrete versions of 
                        A. Schoen's F-Rd surface and I-Wp surfaces.}
\label{fig:example14}
\end{figure}

The fundamental piece $\mathcal T$ of this I-Wp example has six vertices 
\[ p_1=(b,0,b), \;\; p_2=(b,0,0), \;\; p_3=(b,b,0), \]\[ p_4=(1,1,a), \;\; 
p_5=(1,a,1), \;\; p_6=\overline{p_1p_4} \cap \overline{p_3p_5}\; . \]  
There are five triangles in $\mathcal T$, which are 
\[ (p_1,p_2,p_3), \;\; (p_1,p_3,p_6), \;\; (p_3,p_4,p_6), \;\; 
(p_4,p_5,p_6), \;\; (p_5,p_1,p_6). \]  By including the two images of 
$\mathcal T$ under the two reflections across the planes 
$\{ (x,x,z) \, | \, x,z \in R \}$ and 
$\{ (x,y,x) \, | \, x,y \in R \}$, 
we have a larger discrete surface with fifteen triangles.  Reflecting this 
larger piece across all planes of the form 
$\{ (x,y,k) \, | \, x,y \in \mathbb{R} \}$, 
$\{ (x,k,z) \, | \, x,z \in \mathbb{R} \}$, 
$\{ (k,y,z) \, | \, y,z \in \mathbb{R} \}$ for all integers $k$, 
we arrive at a complete embedded triply-periodic surface in $\mathbb{R}^3$.  
See the right-hand side of Figure \ref{fig:example14}.  

The minimality condition that Equation \eqref{eqn:gradientiszero} holds at each 
vertex of the complete triply-periodic surface is 
\begin{equation}\label{whatever3} 1+a+a^2-3 b-2 a b+2 b^2=0 \; , \end{equation} 
\begin{equation}\label{whatever4} a^2+a (2-3 b) + 
b \left( 3b-3+\sqrt{(1+a-b)^2+2(1-b)^2} \right) = 0 \; . \end{equation} 
Thus, to make the surface minimal, we must find $a$ and $b$ satisfying 
Equations \eqref{whatever3}-\eqref{whatever4}, so the minimality 
condition is two-dimensional.  
Equation \eqref{whatever3} holds if 
\[ a= \frac{1}{2} \left( 2 b-1+\sqrt{-3+8 b-4 b^2} \right) \; , \] and then 
Equation \eqref{whatever4} will hold if $b$ satisfies 
\[ \left( 3-\sqrt{-3+8 b-4 b^2} \right) (1-b) = 
\sqrt{2} \sqrt{3-4 b+2 b^2+\sqrt{-3+8 b-4 b^2}} \; . \]  
One can find such a real number $b$ in a completely explicit form (although 
not in such simple forms like in Subsections \ref{section3-3}, 
\ref{section4-2} and \ref{section5-1}).  

One can similarly find a discrete analog, shown on the left-hand side of Figure 
\ref{fig:example14}, of the smooth triply-periodic 
minimal F-Rd surface.  With the simplicial structure chosen in Figure 
\ref{fig:example14}, one can again explicitly 
solve the minimality condition, in the same way as we did for the I-Wp example.  

\subsection{Trigonal example}\label{section5-3} 

\begin{figure}[tbp]
\centerline{
        \hbox{
		\psfig{figure=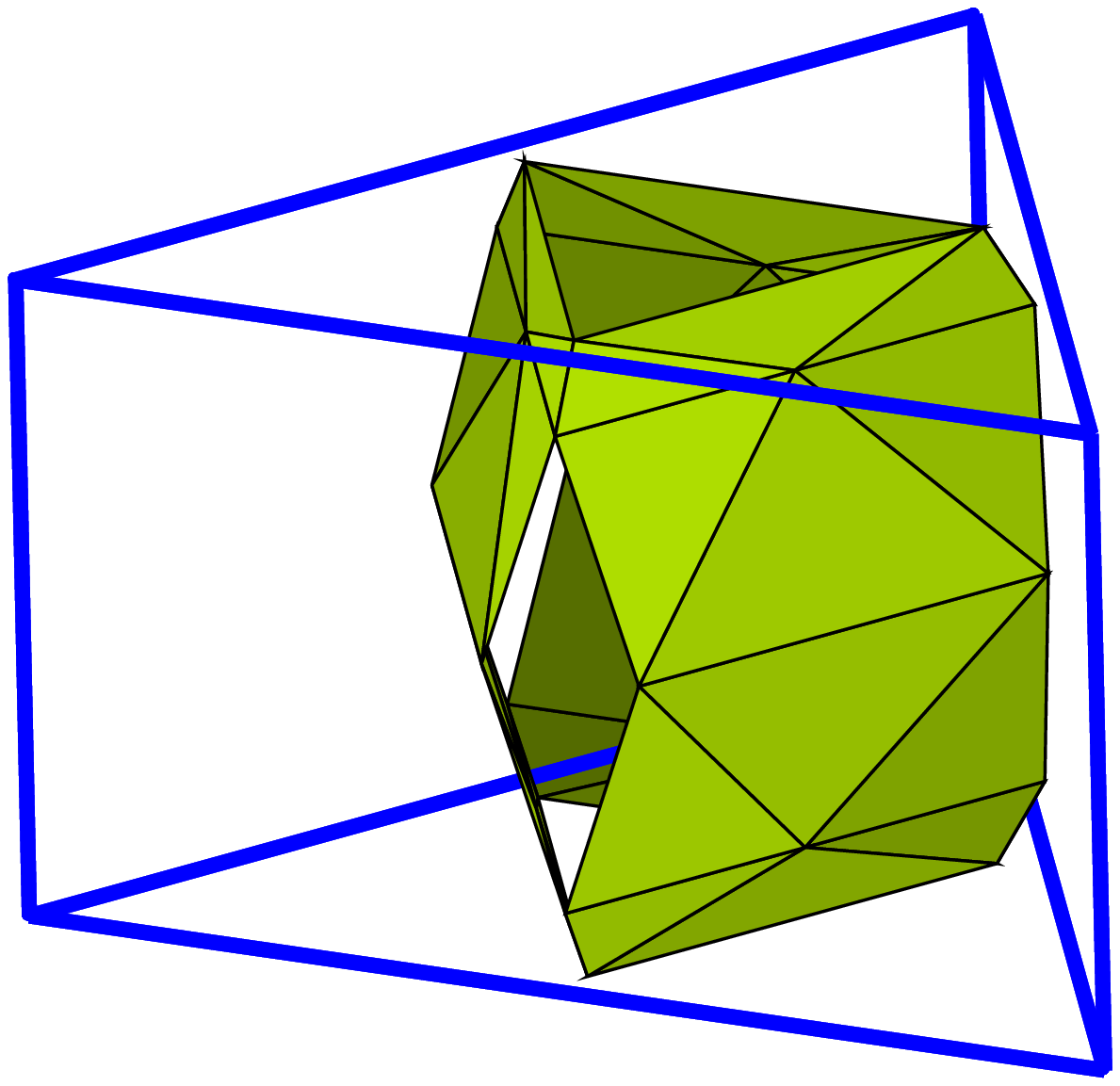,width=1.7in}
		\psfig{figure=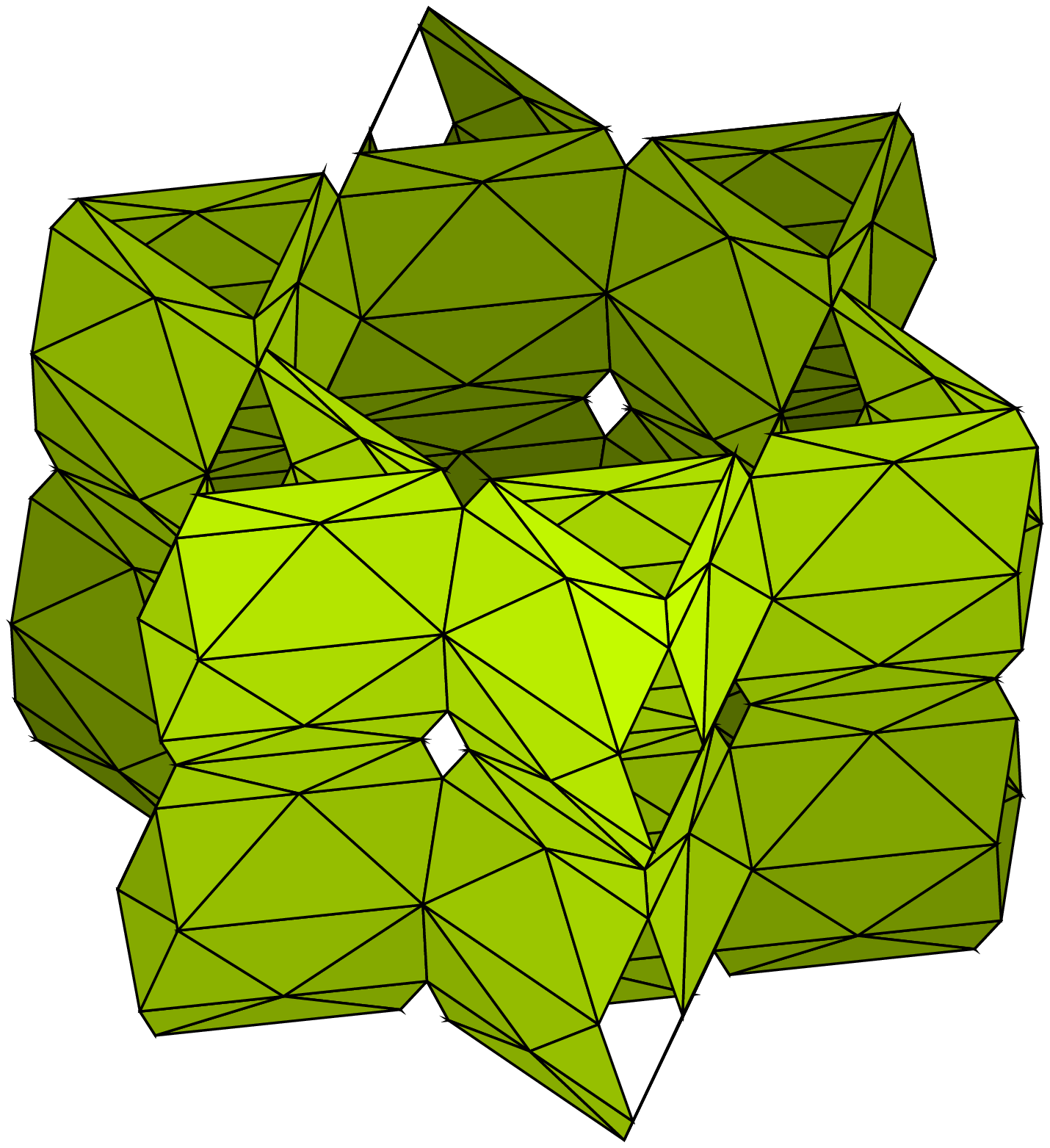,width=1.7in}
		\psfig{figure=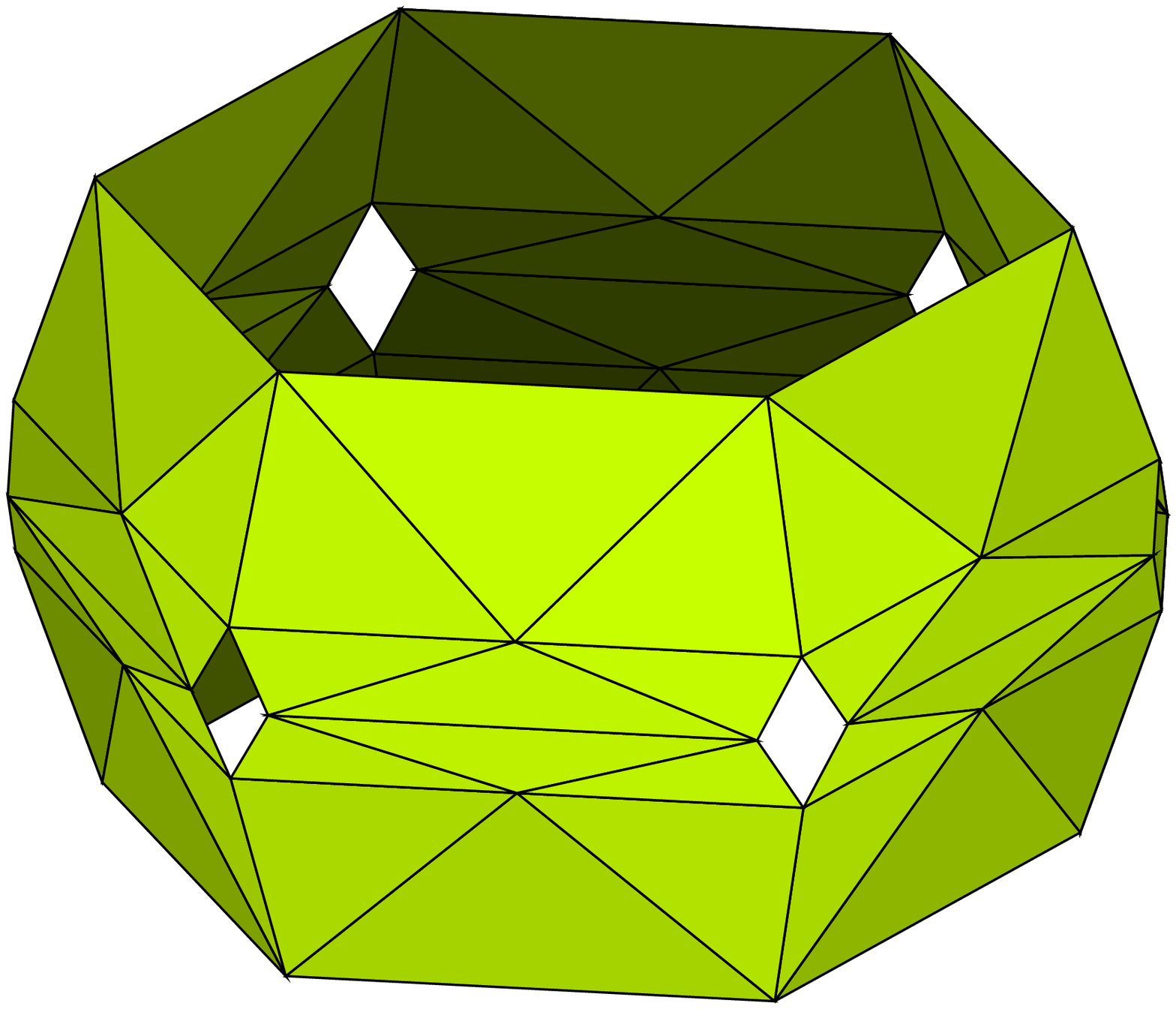,width=1.7in}
	}
}
\caption{\protect\small Discrete version of 
                        A. Schoen's H-T surface.}
\label{fig:example18}
\end{figure}

The fundamental piece $\mathcal T$ of this H-T example has six boundary vertices 
\[ p_1=(\frac{a}{2},\frac{\sqrt{3} a}{2},b), \;\; 
p_2=(\frac{1}{2},\frac{\sqrt{3}}{2},c), \;\; 
p_3=(2-\frac{3 s}{2},\frac{\sqrt{3} s}{2},0), \]
\[ p_4=(2-\frac{3 s}{2},-\frac{\sqrt{3} s}{2},0), \;\; 
p_5=(\frac{1}{2},-\frac{\sqrt{3}}{2},c), \;\; 
p_6=(\frac{a}{2},-\frac{\sqrt{3} a}{2},b). \]
for any given fixed $b>0$, and has one interior vertex 
\[ p_7=\frac{1}{2} (p_2+p_5) \; . \]  There are 
six triangles in $\mathcal T$, which are 
\[ (p_j,p_{j+1},p_7), \; j=1,...,5 \; , \;\;\; (p_6,p_1,p_7). \]
Including the images of $\mathcal T$ under the $120^\circ$ and $240^\circ$ 
rotations about the axis 
$\{ (0,0,z) \, | \, z \in \mathbb{R} \}$, 
and also including the images of $\mathcal T$ and these two rotated copies 
of $\mathcal T$ under reflection across the plane 
$\{ (x,y,0) \, | \, x,y \in \mathbb{R} \}$, 
one has the larger piece shown on the left-hand side of 
Figure \ref{fig:example18}.  This larger piece has five boundary 
components, each contained in a plane, and these five planes bound a 
trigonal prism (a prism of height $2 b$ over an equilateral triangle with 
edge-lengths $2 \sqrt{3}$).  Including the 
images of this larger piece by reflecting across these five planes, and 
also by including all subsequent images of reflections across planes 
containing subsequent boundary components, one arrives at a triply-periodic 
discrete surface, which is embedded when $a,s \in (0,1)$ and $c \in (0,b)$.  
A larger portion of this complete discrete surface is shown in the 
central figure of Figure \ref{fig:example18}.  

The minimality condition involves three equations in the three variables 
$a,c,s$, and so is three-dimensional.  We will not show the 
equations here, but they can be solved 
explicitly.  For example, when 
$b=1$, the following choices ensure minimality: 
\[ a=s=\frac{2+\sqrt{2}}{4} \; , \;\;\; c=\frac{3}{4} \;\;\;\;\;\; (\text{and} 
\;\; b=1) \; . \]  
In fact, these choices also ensure that all of the vertices of $\mathcal T$ 
lie in the same plane, and hence the fundamental piece $\mathcal T$ is planar 
and could be freely triangulated within its trace (see Remark 
\ref{flatpartsofminimalsurfaces}).  

Furthermore, rather than using a portion of the complete surface within a 
trigonal prism as a building block for the complete surface, one could have 
instead used a portion within a hexagonal prism as the building block.  
The figure on the 
right-hand side of Figure \ref{fig:example18} will produce exactly the 
same complete surface (again by reflecting 
across planes containing boundary components).  In the case of 
smooth H-T surfaces, this same duality exists between building blocks in 
trigonal and hexagonal prisms, as noted in \cite{K3}.  

\subsection{Discrete Fischer-Koch example\label{section5-4}} 

\begin{figure}[tbp]
\centerline{
        \hbox{
		\psfig{figure=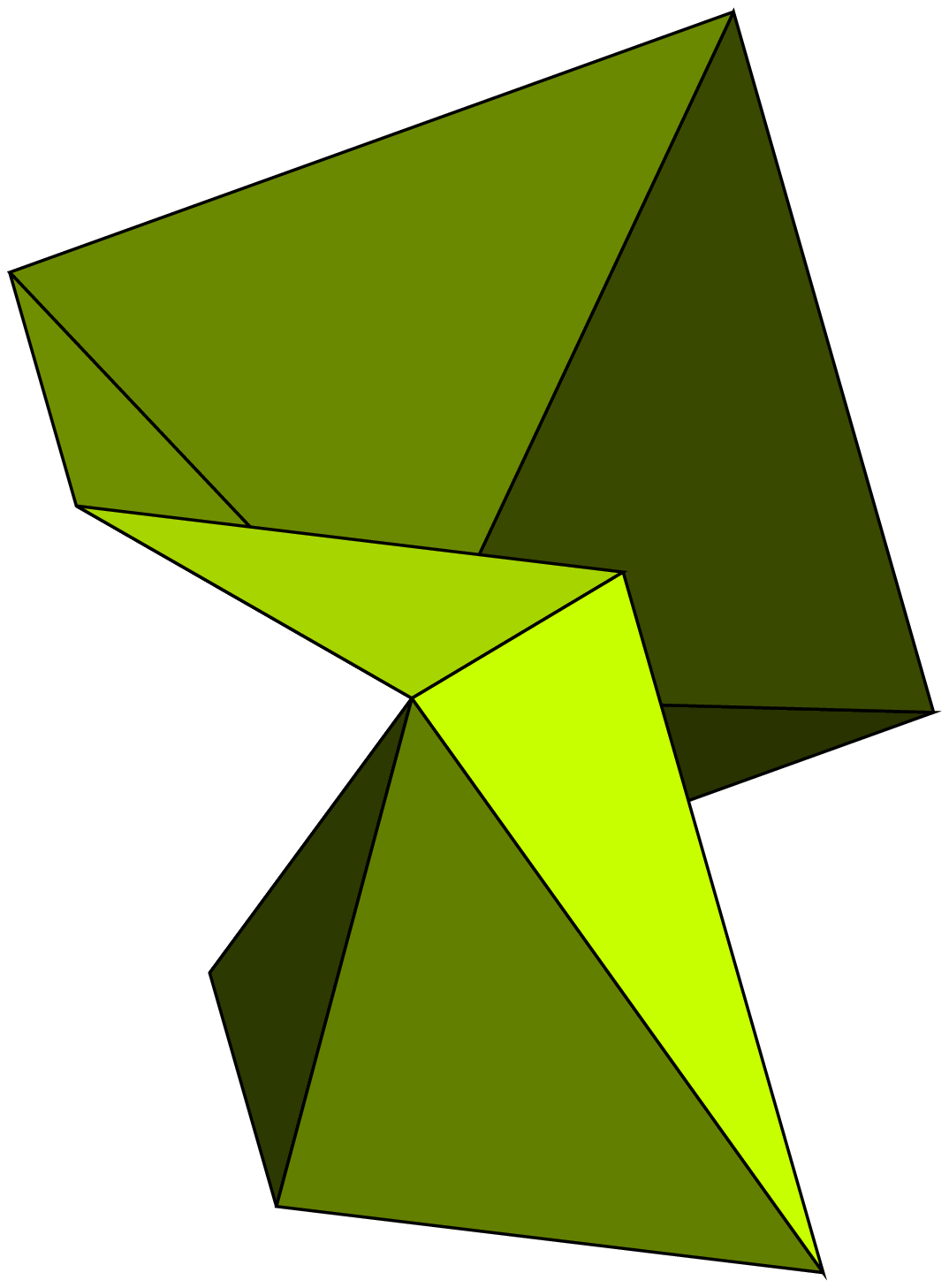,width=1.7in}
		\psfig{figure=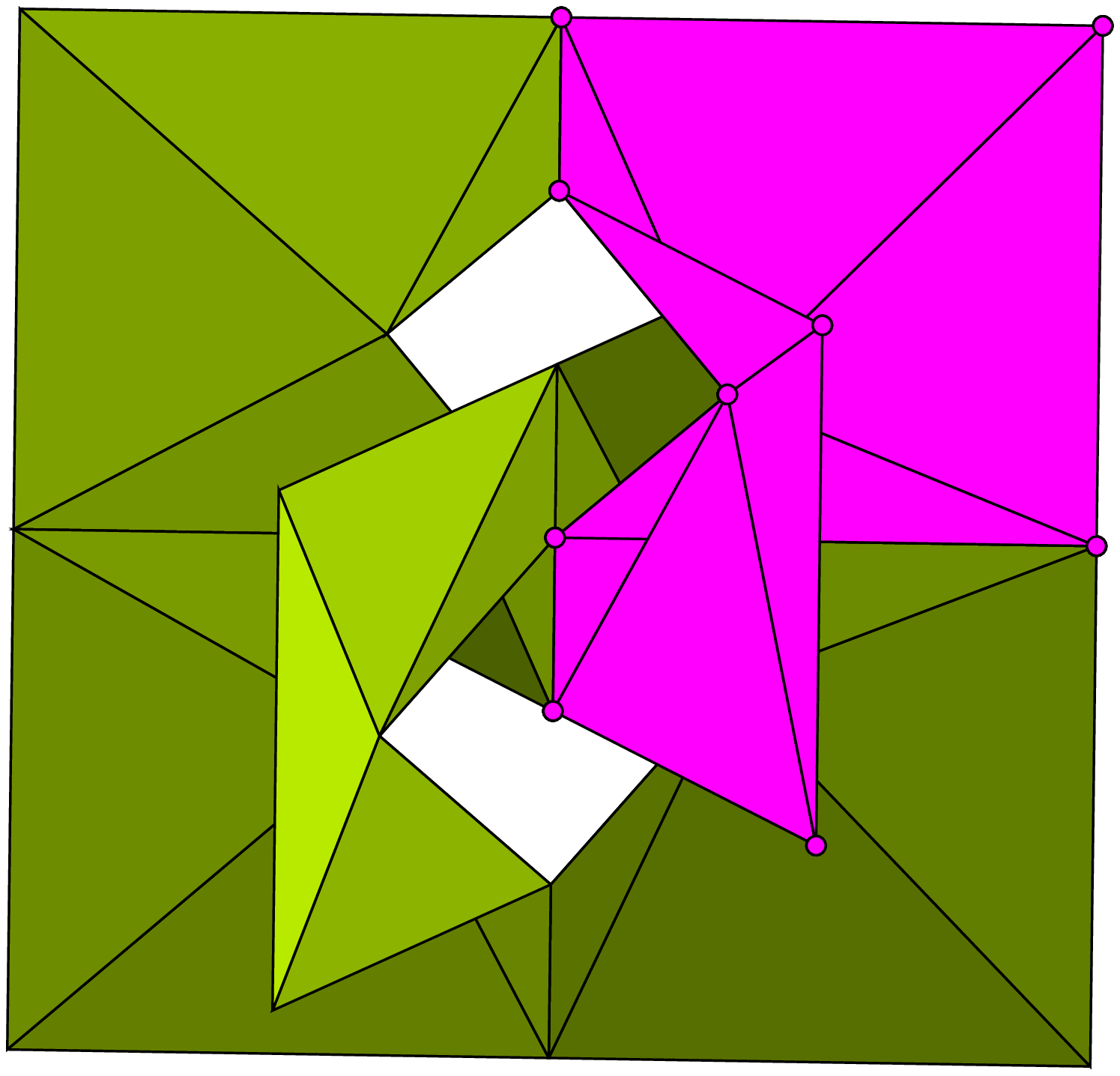,width=1.7in}
		\psfig{figure=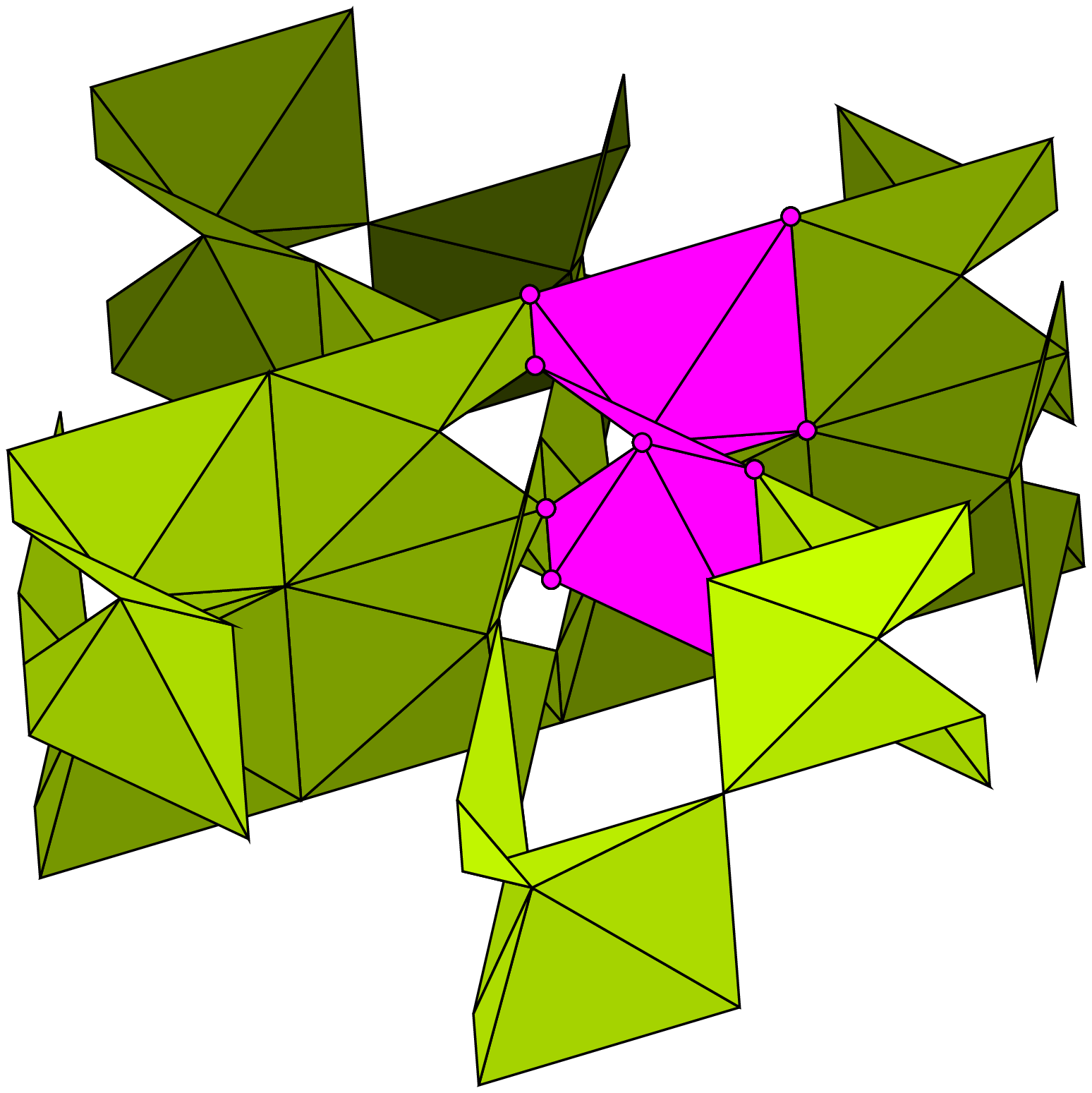,width=1.7in}
	}
}
\centerline{
        \hbox{
		\psfig{figure=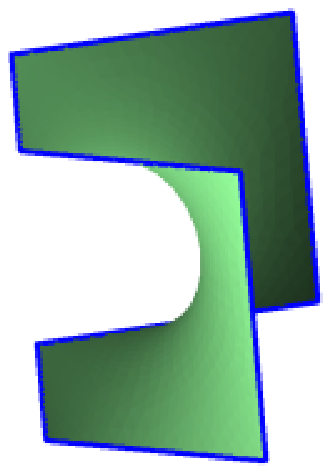,width=1.7in}
		\psfig{figure=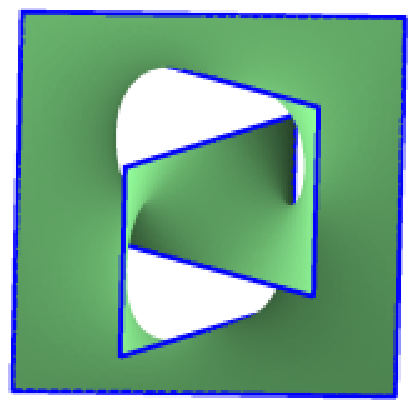,width=2.1in}
	}
}
\caption{\protect\small Discrete and smooth 
                        triply-periodic Fischer-Koch surfaces.}
\label{fig:example10}
\end{figure}

An interesting triply-periodic smooth embedded 
minimal surface was found recently by W. Fischer 
and E. Koch \cite{FC}, and is shown in 
the bottom row of Figure \ref{fig:example10}.  Here we give 
a discrete minimal analog of this surface, shown in 
the top row of Figure \ref{fig:example10}.  

The fundamental piece $\mathcal T$ of this example has eight boundary vertices 
\[ p_1=(0,0,-1), \;\; p_2=(0,0,-2), \;\; p_3=(a,0,-2), \;\; p_4=(a,0,1), \]
\[ p_5=(0,0,1), \;\; p_6=(0,0,2), \;\; p_7=(\frac{a}{2},\frac{\sqrt{3} a}{2},2), 
\;\; p_8=(\frac{a}{2},\frac{\sqrt{3} a}{2},-1). \]
for any given fixed $a>0$, and has one interior vertex 
\[ p_9=(\frac{\sqrt{3} b}{2},\frac{b}{2},0) \; \] with $0<b<a$.  There are 
eight triangles in $\mathcal T$, which are 
\[ (p_j,p_{j+1},p_9), \; j=1,...,7 \; , \;\;\; (p_8,p_1,p_9). \]
The complete triply-periodic surface is 
generated by $180^\circ$ rotations about boundary edges.  
In this example, the symmetry 
Equation \eqref{eqn:gradientiszero} holds at each vertex $p_1$,...,$p_8$ 
in the resulting complete triply-periodic surface.  
However, getting this to hold at $p_9$ requires the proper choice of 
$b$. This minimality condition at $p_9$ is one-dimensional, and one can 
prove existence of a value $b \in (0,a)$ solving it.  

\bibliographystyle{abbrv}

\end{document}